\documentclass[11pt]{amsart}

\usepackage{amsmath,amsfonts,amssymb,amsthm}
\usepackage{pinlabel}
\usepackage{graphicx}
\usepackage[all]{xy}
\usepackage{hyperref}
\xyoption{dvips}

\title{The cylindrical contact homology of universally tight sutured contact solid tori}

\author{Roman Golovko}

\address{Universit\'{e} Paris-Sud, D\'{e}partement de Math\'{e}matiques, Bat. 425, 91405
Orsay, France} \email{roman.golovko@math.u-psud.fr}



\newtheorem{theorem}{Theorem}[section]
\newtheorem{lemma}[theorem]{Lemma}
\newtheorem{fact}[theorem]{Fact}
\newtheorem{claim}[theorem]{Claim}

\newtheorem{proposition}[theorem]{Proposition}
\theoremstyle{remark}
\newtheorem{remark}[theorem]{Remark}
\theoremstyle{definition}
\newtheorem{definition}[theorem]{Definition}

\numberwithin{equation}{subsection}

\begin{document}

\begin{abstract}
We calculate the sutured version of cylindrical contact homology of
a sutured contact solid torus $(S^1\times D^2,\Gamma, \xi)$, where
$\Gamma$ consists of $2n$ parallel sutures of arbitrary slope and
$\xi$ is a universally tight contact structure. In particular, we
show that it is non-zero. This computation is one of the first
computations of the sutured version of cylindrical contact homology
and does not follow from computations in the closed case.
\end{abstract}

\maketitle

\section{Introduction}
The cylindrical contact homology of a (closed) contact manifold was
introduced by Eliashberg and Hofer and is the simplest version of
the symplectic field theory of Eliashberg, Givental and
Hofer~\cite{EliashbergGiventalHofer}. It is the homology of a
differential graded module whose differential counts genus zero
holomorphic curves in the symplectization with one positive puncture
and one negative puncture.

Gabai in the early 1980's
developed the theory of sutured manifolds, see \cite{Gabai}. It became
a powerful tool in studying $3$-manifolds with boundary.
It turns out that there is a way to generalize cylindrical contact
homology to sutured manifolds. This is possible by imposing a
certain convexity condition on the contact form. This construction
is described in the paper of Colin, Ghiggini, Honda and Hutchings
\cite{ColinGhigginiHondaHutchings} and will be summarized in
Section~\ref{section:background}.

In this paper, we construct a sutured contact solid torus with $2n$
parallel sutures of slope $-\frac{k}{l}$ using the gluing method of
Colin, Ghiggini, Honda and Hutchings
\cite{ColinGhigginiHondaHutchings}, and calculate the sutured
cylindrical contact homology of it. Here $n\in \mathbb N$, $(k,l)=1$
and $k>l>0$. In order to define the slope, we choose an oriented
identification $\partial(S^1\times D^2)\simeq T^2 = (\mathbb R /
\mathbb Z)^2$ as follows: map $\{ pt \} \times
\partial D^2$ (the meridian) to $(1, 0)$ (slope is $0$) and $S^1\times \{ pt \}$ (a longitude) to
$(0, 1)$.

This calculation, together with the calculation of the sutured
cylindrical contact homology of the sutured contact solid torus with
$2n$ parallel longitudinal sutures, where $n\geq 2$, that has been
done in~\cite{Golovko}, finishes the calculation of the cylindrical
contact homology of $(S^1\times
D^2,\Gamma, \xi)$, where $\Gamma$ consists of
$2n$
parallel sutures of arbitrary slope, $\xi$ is a universally tight contact structure and such that if one
cuts along the meridian disk, the sutures on the disk are $\partial$-parallel.
In particular, this gives a complete calculation of the cylindrical
contact homology of $(S^1\times
D^2,\Gamma, \xi)$, where $\Gamma$ consists of
$2$
parallel sutures of arbitrary slope and $\xi$ is a universally tight contact structure (observe that in this situation there are only two isomorphic (but not isotopic) universally tight contact structures, see Section 2 in \cite{Honda}). These are not all the universally tight contact structures on the solid torus, but all of them can be obtained from the $\#\Gamma=2$ case by successively applying the folding operation.

Our goal is to prove the following theorem:
\begin{theorem}\label{mainresult}
Let $(S^1\times D^2, \Gamma)$ be a sutured manifold, where $\Gamma$
is a set of $2n$ parallel closed curves of slope $-\frac{k}{l}$,
where $(k,l)=1$, $k>l>0$ and $n\in \mathbb N$. Then there is a
contact form $\alpha$ which makes $(S^1\times D^2, \Gamma, \alpha)$
a sutured contact manifold with a universally tight contact
structrure $\xi=ker\,\alpha$, $HC^{cyl}(S^{1}\times
D^2,\Gamma,\alpha)$ is defined, is independent of the contact form
$\alpha$ for $\xi=ker\,\alpha$ and the almost complex structure $J$
and
\begin{align*}
HC^{cyl,h}(S^{1}\times D^2,\Gamma,\xi)\simeq\left \{
\begin{array}{ll}
\mathbb Q, & \mbox{for} \ k \nmid h > 0;\\
\mathbb Q^{n-1}, & \mbox{for}\ k \mid h > 0;\\
0, & \mbox{otherwise}.
\end{array}
\right.
\end{align*} Here $h$ corresponds to
the homological grading.
\end{theorem}

\section{Background}
\label{section:background}

The goal of this section is to review definitions of sutured contact
manifold and the relative version of cylindrical contact homology.
This section can be considered as a summary of
\cite{ColinGhigginiHondaHutchings}.

\subsection{Review of sutured contact manifolds}
\label{section:glsutcontman}

In this section, we recall some definitions and describe some constructions
from~\cite{ColinGhigginiHondaHutchings}. We first start with the notion of a Liouville manifold.
\begin{definition}
A {\em Liouville manifold} (sometimes called a Liouville domain) is
a pair $(W, \beta)$ which consists of a compact, oriented
$2n$-dimensional manifold $W$ with boundary and a $1$-form $\beta$
on $W$, where $\omega = d\beta$ is a positive symplectic form on $W$
and the {\em Liouville vector field} $Y$ given by $i_{Y}(\omega) =
\beta$ is positively transverse to $\partial W$. It follows that the
$1$-form $\beta_{0} = \beta|_{\partial W}$ is a positive contact form with
kernel $\zeta$.
\end{definition}
We now recall the definition of a sutured contact manifold.
\begin{definition}
A compact oriented $2n+1$-dimensional manifold $M$ with boundary and
corners is a {\em sutured contact manifold} if it comes with an oriented,
not necessarily connected submanifold $\Gamma \subset \partial M$ of
dimension $2n-1$, called the {\em suture}, together with a
neighborhood $U(\Gamma)=[-1,0]\times[-1, 1]\times \Gamma$ of $\Gamma
= \{0\}\times \{0\}\times \Gamma$ in $M$, with coordinates
$(\tau,t)\in [-1,0]\times [-1,1]$, such that the following holds:
\begin{itemize}
\item[($1$)] $U\cap \partial M=(\{0\}\times [-1,1]\times \Gamma)\cup([-1,0]\times
\{-1\} \times \Gamma)\cup ([-1,0]\times \{1\}\times \Gamma)$;
\item[($2$)] $\partial M \setminus (\{0\}\times (-1,1)\times \Gamma)=R_{-}(\Gamma)\sqcup
R_{+}(\Gamma)$, where the orientation of $\partial M$ coincides with
that of $R_{+}(\Gamma)$ and is opposite that of $R_{-}(\Gamma)$ and
the orientation of $\Gamma$ coincides with the boundary orientation of
$R_{+}(\Gamma)$;
\item[($3$)] the corners of $M$ are precisely $\{0\}\times\{\pm 1\}\times
\Gamma$.
\end{itemize}
In addition, $M$ is equipped with a
contact structure $\xi$, which is the kernel of a positive
contact 1-form $\alpha$ such that:
\begin{itemize}
\item[($i$)] $(R_{\pm}(\Gamma), \beta_{\pm}=\alpha|_{R_{\pm}(\Gamma)})$
is a
Liouville manifold;

\item[($ii$)] $\alpha=Cdt+\beta$ inside $U(\Gamma)$, where $C>0$ and $\beta$ is independent of
$t$ and does not have a $dt$-term;

\item[($iii$)] $\partial_{\tau} = Y_{\pm}$, where $Y_{\pm}$ is a Liouville vector field for $\beta_{\pm}$.
\end{itemize}
Such a contact form $\alpha$ is called {\em adapted} to
$(M,\Gamma,U(\Gamma))$.
\end{definition}

Here we briefly describe the way to glue sutured contact
manifolds. This procedure was
first described by Colin and Honda in \cite{ColinHonda} and then generalized by Colin, Ghiggini, Honda and Hutchings
in~\cite{ColinGhigginiHondaHutchings}.

Let $(M',\Gamma',U(\Gamma'),\xi')$ be a sutured contact $3$-manifold
with an adapted contact form $\alpha'$. Let $\pi$ be the
projection along $\partial_{t}$ defined on $U(\Gamma')$.

Consider $2$-dimensional submanifolds $P_{\pm}\subset R_{\pm}(\Gamma')$
such that $\partial P_{\pm}$ is the union of $(\partial
P_{\pm})_{\partial}\subset \partial R_{\pm}(\Gamma')$, $(\partial
P_{\pm})_{int}\subset int(R_{\pm}(\Gamma'))$ and $\partial P_{\pm}$
is positively transversal to the Liouville vector field $Y'_{\pm}$
on $R_{\pm}(\Gamma')$. When we write $(\partial P_{\pm})_{int}$ and $(\partial
P_{\pm})_{\partial}$, we assume that closures are taken as
appropriate. In addition, we assume that $\pi((\partial
P_{-})_{\partial})\cap \pi(\partial P_{+})_{\partial})=\emptyset$.

Consider a diffeomorphism $\varphi$ which sends
$(P_{+},\beta'_{+}|_{P{+}})$ to $(P_{-},\beta'_{-}|_{P_{-}})$, $(\partial P_{+})_{int}$ to $(\partial P_{-})_{\partial}$ and
$(\partial P_{+})_{\partial}$ to $(\partial P_{-})_{int}$. Note that, since dim $M = 3$, we
only need $\beta'_{+}|_{P_{+}}$ and
$\varphi^{\ast}(\beta'_{-}|_{P_{-}})$ to match up on $\partial
P_{+}$, since we can linearly interpolate between primitives of
positive area forms on a surface.

Topologically, $(M, \Gamma)$ is constructed from
$(M', \Gamma')$ and the gluing data $(P_{+}, P_{-}, \varphi)$ as
follows: Let $M=M'/\sim$, where
\begin{itemize}
\item $x\sim \varphi(x)$ for all $x\in P_{+}$;

\item $x\sim x'$ if $x, x'\in \pi^{-1}(\Gamma')$ and $\pi(x)=\pi(x')\in
\Gamma'$.
\end{itemize}
Then
\begin{align*}
R_{\pm}(\Gamma)=\frac{\overline{R_{\pm}(\Gamma')\setminus
P_{\pm}}}{(\partial P_{\pm})_{int}}\sim \pi_{\pm}((\partial
P_{\mp})_{\partial})
\end{align*}
and
\begin{align*}
\Gamma=\frac{\overline{\Gamma'\setminus \pi (\partial P_{+}\sqcup
\partial P_{-})}}{\pi((\partial P_{+})_{int}\cap (\partial
P_{+})_{\partial})} \sim \pi((\partial P_{-})_{int}\cap (\partial
P_{-})_{\partial}).
\end{align*}

For the detailed description of the gluing procedure we refer to
\cite{ColinGhigginiHondaHutchings}.

Finally, we describe the way to complete sutured contact manifold $(M, \alpha)$
to a noncompact contact manifold $(M^{\ast}, \alpha^{\ast})$. This construction was first described in~\cite{ColinGhigginiHondaHutchings}.

Given a sutured contact manifold $(M,\Gamma, U(\Gamma),\xi)$ with
an adapted contact form $\alpha$. The form $\alpha$ is then defined by
$Cdt +\beta_{\pm}$ on $[1-\varepsilon,1]\times
R_{+}(\Gamma)$ and $[-1,-1 +\varepsilon]\times R_{-}(\Gamma)$ of
$R_{+}(\Gamma)=\{1\}\times R_{+}(\Gamma)$ and $R_{-}(\Gamma) =
\{-1\} \times R_{-}(\Gamma)$, where $t\in [-1,-1+\varepsilon]\cup
[1-\varepsilon,1]$ extends the $t$-coordinate on $U$. On $U$,
$\alpha$ is given by $Cdt + \beta$, $\beta = \beta_{+} = \beta_{-}$ and
$\partial_{\tau}$ is a Liouville vector field $Y$ for $\beta$.
We now extend $\alpha$ to $[1,\infty)\times R_{+}(\Gamma)$ and
$(-\infty,-1]\times R_{-}(\Gamma)$ by taking $Cdt + \beta_{\pm}$ as
appropriate. The boundary of this new manifold is $\{0\}\times
\mathbb R\times \Gamma$. Since $\partial_{\tau} = Y$,
the form $d\beta|_{[-1,0]\times \{t\}\times \Gamma}$ coincides with the
symplectization of $\beta|_{\{0\}\times \{t\}\times \Gamma}$ in the
positive $\tau$-direction. We then glue $[0,\infty)\times \mathbb R
\times \Gamma$ with the form $Cdt+e^{\tau}\beta_{0}$, where
$\beta_{0}$ is the pullback of $\beta$ to $\{0\}\times \{t\}\times
\Gamma$.

We denote by $M^{\ast}$ the noncompact extension of $M$ (described above)
and by $\alpha^{\ast}$ the extension of $\alpha$ to $M^{\ast}$.

\subsection{Review of cylindrical contact homology}
\label{section:sutcylconthom}
In this section, we review the definition of cylindrical contact homology for sutured manifolds. We refer to \cite{ColinGhigginiHondaHutchings} for more details of this construction.

Let $(M, \Gamma, U(\Gamma), \xi)$ be a sutured contact manifold with
an adapted contact form $\alpha$ and $(M^{\ast}, \alpha^{\ast})$ be its
completion.

The {\em Reeb vector field} $R_{\alpha^{\ast}}$ associated
to $\alpha^{\ast}$ is given by $d\alpha^{\ast}(R_{\alpha^{\ast}},\cdot)=0$ and
$\alpha^{\ast}(R_{\alpha^{\ast}})=1$.
We assume that
$R_{\alpha^{\ast}}$ is {\em nondegenerate}, i.e., the first return map along each (not necessarily simple)
periodic orbit does not have 1 as an eigenvalue. Observe that nondegeneracy can always be achieved by a small
perturbation.

\begin{remark}\label{reeborbincompl}
Note that every periodic orbit of $R_{\alpha^{\ast}}$ lies in $M$. Thus, the
set of periodic Reeb orbits of $R_{\alpha^{ \ast}}$ coincides with
the set of periodic Reeb orbits of $R_{\alpha}$.
\end{remark}

A Reeb orbit $\gamma$ is called {\em elliptic} or {\em positive}
({\em negative}) {\em hyperbolic} if the eigenvalues of
$P_{\gamma}$ lie on the unit circle or the positive (negative)
real line.

If $\tau$ is a trivialization of $\xi$ over $\gamma$, we can then
define the Conley-Zehnder index. In $3$-dimensional case,
we can explicitly describe the Conley-Zehnder
index and its behavior under multiple covers in the following way:

\begin{proposition}[\cite{Hutchings}]
If $\gamma$ is an elliptic orbit, then there exists an irrational number $\phi\in
\mathbb R$ such that $P_{\gamma}$ is conjugate in $SL_2(\mathbb R)$
to a rotation by angle $2\pi\phi$ and
\begin{align*}
\mu_{\tau}(\gamma^k)=2\lfloor k\phi \rfloor + 1,
\end{align*}
where $2\pi\phi$ is the total rotation angle with respect to $\tau$
of the linearized flow around the orbit.

If $\gamma$ is a positive (negative) hyperbolic, then
there is an even (odd) integer $r$ such that the
linearized flow around the orbit rotates the eigenspaces of
$P_{\gamma}$ by angle $\pi r$ with respect to $\tau$ and in this case
\begin{align*}
\mu_{\tau}(\gamma^k)=kr.
\end{align*}
\end{proposition}

A closed orbit of $R_{\alpha^{ \ast}}$ is said to be {\em
good} if it does not cover a simple orbit $\gamma$ an even number of times,
where the first return map $\xi_{\gamma(0)}\to \xi_{\gamma(T)}$ has an odd number of eigenvalues in the interval $(-1,0)$.
Here $T$ is the period of the orbit $\gamma$.
An orbit that is not good is called {\em
bad}.

We now recall the notion of an almost complex structure on $\mathbb R\times M^{\ast}$ that is tailored to $(M^{\ast}, \alpha^{\ast})$.

Let $(W, \beta)$ be a Liouville manifold and $\zeta$ be the contact
structure given on $\partial W$ by $ker(\beta_{0})$, where
$\beta_{0} = \beta|_{\partial W}$. Also, let $(\widehat{W},
\widehat{\beta})$ be the completion of $(W, \beta)$, namely
$\widehat{W} = W\cup ([0,\infty)\times \partial W)$ and $\widehat
{\beta}|_{[0,\infty)\times \partial W} = e^{\tau}\beta_{0}$, where
$\tau$ is the $[0,\infty)$-coordinate. An almost complex structure
$J_{0}$ on $\widehat W$ is {\em $\widehat{\beta}$- adapted} if
$J_{0}$ is $\beta_{0}$-adapted on $[0,\infty)\times \partial W$; and
$d\beta (v,J_{0}v) > 0$ for all nonzero tangent vectors $v$ on $W$.

\begin{definition}
Given a sutured contact manifold $(M, \Gamma, U(\Gamma), \xi)$. Let
$\alpha$ be an adapted contact form and $(M^{\ast}, \alpha^{\ast})$
be its completion. We say that an almost complex structure $J$ on
$\mathbb R\times M^{\ast}$ is {\em tailored to} $(M^{\ast},
\alpha^{\ast})$ if the following conditions hold:
\begin{itemize}
\item[($1$)] $J$ is $\alpha^{\ast}$-adapted, i.e., $J$ is $\mathbb R$-invariant, $J(\xi)=\xi$, $d\alpha(v,Jv) > 0$ for nonzero $v\in \xi$ and $J(\partial_s)=R_{\alpha^{\ast}}$, where $s$ is the $\mathbb R$-coordinate;
\item[($2$)] $J$ is $\partial_t$-invariant in a neighborhood of $M^{\ast} \setminus int(M)$;
\item[($3$)] The projection of $J$ to $T \widehat{R_{\pm}(\Gamma)}$ is a $\widehat{\beta}_{\pm}$-adapted almost complex
structure $J_{0}$ on the completion $(\widehat{R_{+}(\Gamma)},
\widehat{\beta}_{+}) \bigsqcup (\widehat{R_{-}(\Gamma)},
\widehat{\beta}_{-})$ of the Liouville manifold $(R_{+}(\Gamma),
\beta_{+}) \bigsqcup (R_{-}(\Gamma), \beta_{-})$. Moreover, the flow
of $\partial_t$ identifies $J_0|_{\widehat{R_{+}(\Gamma)}\setminus
R_{+}(\Gamma)}$ and $J_{0}|_{\widehat{R_{-}(\Gamma)}\setminus
R_{-}(\Gamma)}$.
\end{itemize}
\end{definition}

Given a sutured contact manifold $(M, \Gamma, U(\Gamma), \alpha)$ and an $\alpha^{\ast}$-adapted almost complex structure $J$, we
define the {\em sutured cylindrical contact homology group}
$HC^{cyl}(M, \Gamma, \alpha, J)$ to be the cylindrical contact
homology of $(M^{\ast}, \alpha^{\ast}, J)$. The cylindrical contact
homology chain complex $C(\alpha, J)$ is a $\mathbb Q$-module freely
generated by all good Reeb orbits, where the grading $|\cdot|$ and the
boundary map $\partial$ are defined as in~\cite{Bourgeois} with
respect to the $\alpha^{\ast}$-adapted almost complex structure $J$.
The homology of $C(\alpha, J)$ is the sutured cylindrical contact
homology group $HC^{cyl}(M,\Gamma, \alpha, J)$.

For our calculations we need the following construction of a
``global'' symplectic trivialization described in~\cite{Bourgeois}.
Assume that all the Reeb orbits of $R_{\alpha}$ are good. Let us now
choose trivializations $\tau(\gamma)$ consistently for all Reeb
orbits $\gamma$. Assume that $H_1(M;\mathbb Z)$ is a free module. We
pick representatives $C_1,\dots,C_s$ in $H_1(M;\mathbb Z)$ for a
basis of $H_1(M;\mathbb Z)$, together with a trivialization of $\xi$
along each representative $C_i$, $i = 1,\dots,s$.  Now for a Reeb
orbit $\gamma$, we distinguish the following cases:
\begin{itemize}
\item[($1$)] $[\gamma]=0 \in H_1(M;\mathbb Z)$. Choose a spanning surface
$S_{\gamma}$ and use it to trivialize $\xi$ along $\gamma$.

\item[($2$)] $0\ne [\gamma]\in H_1(M;\mathbb Z)$. We choose a surface $S_{\gamma}$ realizing a homology
between $\gamma$ and a linear combination of the representatives
$C_i$, $i=1,\dots,s$. We then use $S_{\gamma}$ to extend the chosen
trivializations of $\xi$ along the $C_i$, $i=1,\dots,s$ to $\gamma$.
\end{itemize}

We denote the obtained trivialization by $\tau$.

To a $J$-holomorphic cylinder in $\mathcal M^{J}(\gamma; \gamma')$,
we can glue the chosen surfaces $S_{\gamma}$ and $S_{\gamma'}$ and
obtain a closed surface in $M$ (here $\mathcal M^{J}(\gamma;\gamma')$ is
a moduli space of $J$-holomorphic cylinders considered in
cylindrical contact homology theory). Let $A\in H_{2}(M;\mathbb Z)$
be its homology class; we can use it to decorate the corresponding
connected component $M^{J}_{A}(\gamma; \gamma')$ of the moduli
space. Using $\tau$ we can write
\begin{align}\label{dimension}
ind(u)=|\gamma|-|\gamma'|+2\langle c_{1}(\xi), A \rangle
\end{align}
for $u\in \mathcal M^{J}_{A}(\gamma; \gamma')$, where $|\gamma|$ is
the {\em Conley-Zehnder grading of $\gamma$} defined by
\begin{align}\label{grading}
|\gamma|:=\mu_{\tau}(\gamma)-1.
\end{align}
We will use Formulas~\ref{dimension} and~\ref{grading} for our
calculations.

In addition, we will need the following fact, which is a consequence
of Lemma 5.4 in \cite{BourgeoisEliashbergHoferWysockiZehnder}:
\begin{fact}\label{energy}
Let $(M,\alpha)$ be a closed, oriented contact manifold with
nondegenerate Reeb orbits and $u\in\mathcal M^{J}(\gamma; \gamma')$,
where $\gamma$ and $\gamma'$ are good Reeb orbits and $J$ is an
$\alpha$-adapted almost complex structure on $\mathbb R\times M$.
Then $\mathcal A(\gamma):= \int_{\gamma}\alpha \geq
\int_{\gamma'}\alpha =: \mathcal A(\gamma')$ with equality if and
only if $\gamma=\gamma'$ and in this case the moduli space consists
of a single element $\mathbb R\times \gamma$.
\end{fact}
Now we recall the following theorem:
\begin{theorem}[\cite{Bourgeois}]\label{welldefconthom}
Let $(M,\alpha)$ be a closed, oriented contact manifold with
nondegenerate Reeb orbits. Let $C_{m}^{h}(M,\alpha)$ be the
cylindrical contact homology complex, where $h$ is a homotopy class
of Reeb orbits and $m$ corresponds to the Conley-Zehnder grading. If
$C_{k}^{0}(M,\alpha)=0$ for $k=-1,0,1$, then for every free homotopy
class $h$
\begin{itemize}

\item[($1$)] $\partial^2=0$;

\item[($2$)] $H(C^{h}_{\ast}(M,\alpha),\partial)$ is independent of the contact form $\alpha$ for $\xi$,
the almost complex structure $J$ and the choice of perturbation for
the moduli spaces.
\end{itemize}
\end{theorem}
When $M$ is closed and $\mathbb R\times M$ is $4$-dimensional, the
following transversality result has been proven by Momin, see Proposition~2.10
in~\cite{Momin}:
\begin{theorem}[\cite{Momin}]\label{transcylinder}
Let $u\in \mathcal M^{J}(\gamma; \gamma')$ be such that $ind(u)=1$.
Then the linearization of the Cauchy-Riemann operator is surjective
at $u$.
\end{theorem}

\begin{remark}
Observe that Theorem~\ref{transcylinder} does not require $J$ to be
generic. In addition, note that Theorem~\ref{transcylinder} can be considered as a
consequence of the automatic transversality result of Wendl, see
Theorem~$0.1$ in~\cite{Wendl}.
\end{remark}

Finally, we recall the following result of Colin, Ghiggini, Honda
and Hutchings from~\cite{ColinGhigginiHondaHutchings}:

\begin{theorem}[\cite{ColinGhigginiHondaHutchings}]\label{fullwelldefconthom}
Let $(M, \Gamma, U(\Gamma), \xi)$ be a sutured contact $3$-manifold
with an adapted contact form $\alpha$, $(M^{\ast}, \alpha^{\ast})$
be its completion and $J$ be an almost complex structure on $\mathbb
R\times M^{\ast}$ which is tailored to $(M^{\ast}, \alpha^{\ast})$.
Then the contact homology algebra $HC(M, \Gamma, \xi)$ is defined
and independent of the choice of contact $1$-form $\alpha$ with $ker
(\alpha) = \xi$, adapted almost complex structure $J$, and abstract
perturbation.
\end{theorem}

\begin{remark}\label{translationforsutures}
Fact~\ref{energy}, Theorems~\ref{welldefconthom}
and~\ref{transcylinder}, Formulas~\ref{dimension} and~\ref{grading}
hold for $J$-holomorphic curves in the symplectization of the
completion of a sutured contact manifold, provided that we choose
the almost complex structure $J$ on $\mathbb R \times M^{\ast}$ to
be tailored to $(M^{\ast}, \alpha^{\ast})$.
\end{remark}

\begin{remark}
Observe that Theorem~\ref{fullwelldefconthom} and
Remark~\ref{translationforsutures} rely on the assumption that the
machinery, needed to prove the analogous properties for contact
homology and cylindrical contact homology in the closed situation, works.
\end{remark}

\section{Construction}\label{constrodthesoltornaint}
The goal of this section is to construct the sutured contact solid
torus $(S^1\times D^2, \tilde{\Gamma}, \tilde{\alpha}_{\delta})$,
where $\tilde{\Gamma}$ consists of $2n$ parallel sutures of slope
$-\frac{k}{l}$, $(k,l)=1$, $k>l>0$ and $n\in \mathbb N$. Here
$\tilde{\alpha}_{\delta}$ is a contact form such that
$\xi=ker\,\tilde{\alpha}_{\delta}$ is a universally tight contact
structure and the set of embedded orbits of
$R_{\tilde{\alpha}_{\delta}}$ consists of an elliptic orbit $\gamma$
and hyperbolic orbits $\gamma_1,\dots,\gamma_n$ with
\begin{align*}
&[\gamma]=1,\ [\gamma_{i}]=k \in \mathbb Z \simeq H_{1}(S^1\times D^2;\mathbb Z),\ \mathcal A(\gamma^{k})>\mathcal A(\gamma_{i}),\\
&\mu_{\tau}(\gamma^{s}_{i})=
-2ls\quad \mbox{and}\\
&\mu_{\tau}(\gamma^t)=-2ml+1,
\end{align*}
where $(m-1)k<t\leq mk$, for  some ``global''  symplectic
trivialization $\tau$. Here
 $i=1,\dots,n$,  $t\leq N_{\delta}$, $s\leq
\frac{N_{\delta}}{k}$, $N_{\delta}\gg 0$.

\subsection{Gluing map} First we construct $H\in C^{\infty}(\mathbb
R^2)$. The time-$1$ flow of the Hamiltonian vector field associated
to $H$ composed with an appropriate rotation will play a role of the
gluing map when we will apply the gluing construction described in
Section~\ref{section:glsutcontman} to the sutured contact solid
cylinder constructed in Section~\ref{section:construction}.

We fix $p\in \mathbb R^2$ and consider $H_{sing}:\mathbb
R^2\rightarrow \mathbb R$ given by $H_{sing}=\mu r^2
\cos(nk\theta)$ in polar coordinates $(r,\theta)$ about $p$, where $\mu>0$,
$n\geq 1$ and $k\in \mathbb N\setminus\{1\}$. Note that
$H_{sing}$ is singular only at $p$.

\begin{lemma}\label{coordinates}
There exists a function $H\in C^{\infty}(\mathbb R^2)$ which satisfies the following properties:
\begin{itemize}
\item $H=H_{sing}$ on $\mathbb R^2\setminus D(r_{sing})$ for some $r_{sing}>0$;
\item $H$ is $\frac{2\pi}{nk}$-symmetric with respect to $\theta$;
\item the set of critical points of $H$ consists of equally spaced saddle points $p_1,\dots,p_{nk}$ and a critical point $p$;
\item  there exists a neighborhood $U_{s}$ of
$p_{s}$ with coordinates $(x,y)$ such that $H=axy$ on $U_s$ with $a>0$, and such that $\frac{2\pi}{nk}$-rotation about $p$ that we call $R_{nk}$ maps $U_s$ with the corresponding coordinate system to $U_{s+1}$ with the corresponding coordinate system for $s=1,\dots,nk$;
\item there exists a neighborhood $U$ of $p$ such that $H=\tilde{B}r^2-\tilde{C}$ on $U$, where $\tilde{C}>0$ and
$\tilde{B}$ is a small positive number.
\end{itemize}
\end{lemma}
\begin{proof}
 We construct $H\in
C^{\infty}(\mathbb R^2)$ from $H_{sing}$ by perturbing $H_{sing}$ on
a disk $D(r_{sing})$ about $p$ in such a way that $H$ has $nk$
equally spaced saddle points, critical point at $p$ and interpolates
with no other critical points with $H_{sing}$. In other
words, $H=H_{sing}$ on $\mathbb R^2\setminus D(r_{sing})$ for some $r_{sing}>0$. For the
level sets of $H_{sing}$ and $H$ in the case $n=1$, $k=3$ we refer
to Figure~\ref{levelsets5}.

The construction of $H$ is a modification of the construction
described in~\cite{Cotton-Clay}.

We proceed in four steps.
\begin{itemize}
\item[($1$)] We consider
\begin{align*}
H_1&=H_{sing}+f(r,\theta)=H_{sing}+f_{\exp}(r,\theta)+g(r,\theta)\\&=\mu r^2\cos(nk\theta)-Ae^{-mr^2}+g(r,\theta),
\end{align*}
where $A$ and $m$ are positive constants, and $g(r,\theta)$ is a smooth
function to be chosen later. We are interested in the critical
points of $H_1$ away from the origin.

\begin{figure}[t]
\includegraphics[width=400pt]{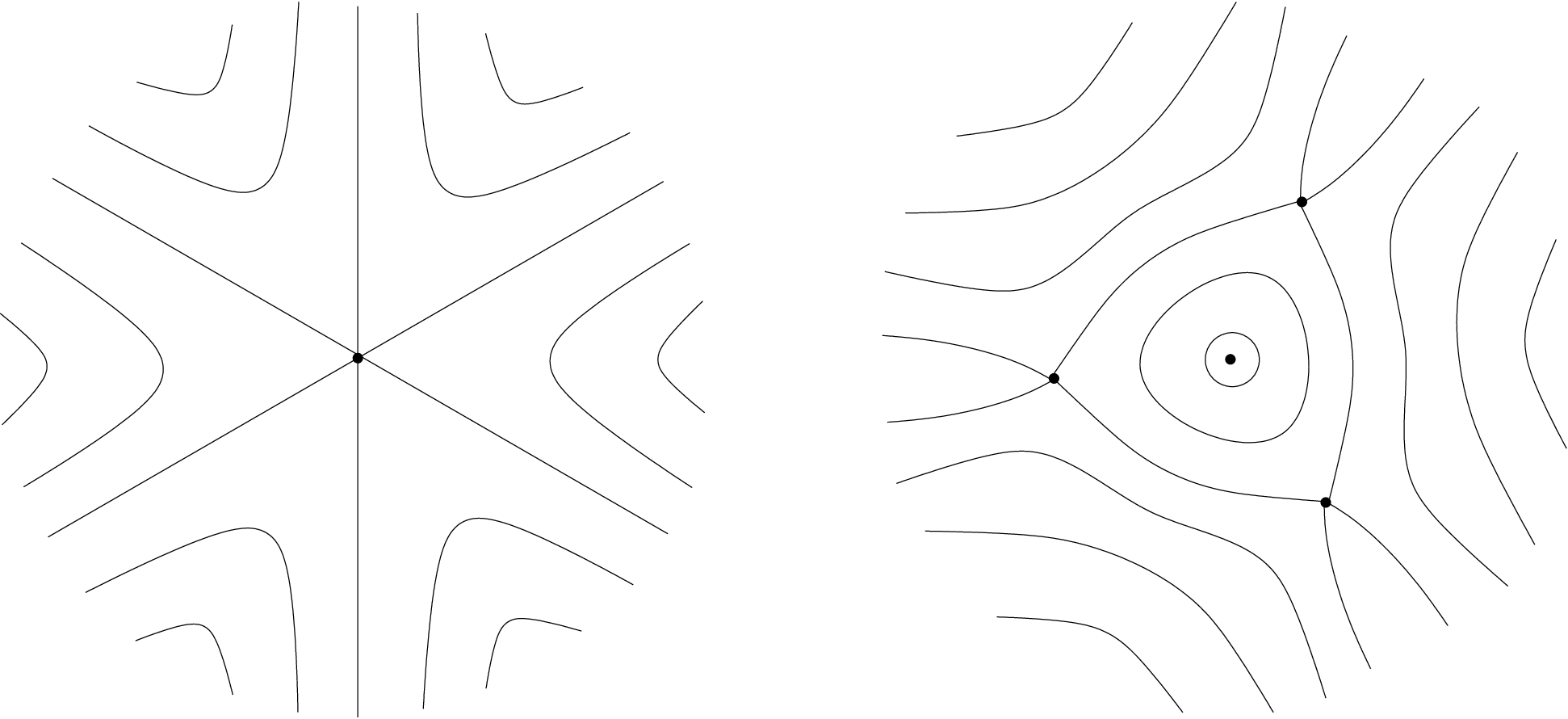}
\caption{The level sets of $H_{sing}$ (left) and the level sets of
$H$ (right) in the case $n=1$, $k=3$} \label{levelsets5}
\end{figure}

We calculate
\begin{align*}
&\frac{\partial H_1}{\partial r}= 2\mu r \cos(nk\theta)+2mrAe^{-mr^2}+\frac{\partial g}{\partial r},\\
&\frac{\partial H_1}{\partial \theta} = -nk\mu r^2
\sin(nk\theta).
\end{align*}
Thus, at the critical points of $H_1$ we must have $\sin(nk\theta)
= 0$. In this case, $\cos(nk\theta) = \pm1$. If $\cos(nk\theta)
= 1$, then $\frac{\partial H_1}{\partial r}-\frac{\partial
g}{\partial r}$ cannot be zero. When $\cos(nk\theta) =-1$,
$\frac{\partial H_1}{\partial r}-\frac{\partial g}{\partial r} =
-2\mu r+2mrAe^{-mr^2}$. For $r
> 0$, $\frac{\partial H_1}{\partial r}-\frac{\partial g}{\partial r}=0$ when $e^{mr^2} =
\frac{mA}{\mu}$, i.e., when $r=r_{c}:=\sqrt{\frac{1}{m}
\ln(\frac{mA}{\mu})}$. We impose the restriction that $mA > \mu$.
Note that by making $m$ large, we can make $r_{c}$ arbitrarily
small. When $\cos(nk\theta)=-1$,
$H_{1}-g(r,\theta)=-\frac{\mu}{m}(\ln(\frac{mA}{\mu})+1)$. Let
$g(r)$ be equal to $\frac{\mu}{m}(\ln(\frac{mA}{\mu})+1)$ on the
annular neighborhood of $r=r_{c}$. For such $g$, $H_{1}$ is $0$ at
the critical points, i.e., at the points $(r_{c},\theta)$, where
$\cos(nk\theta)=-1$.

In summary, we get critical points at one value of $r$ at the values
of $\theta$ when $\cos(nk\theta)=-1$, that is, for $nk$ values
of $\theta$. These are our $nk$ saddle points (it's not hard to
see they are saddle points; alternatively, we can deduce that they
must be for index reasons).

\item[($2$)] Keeping $f_{exp}$ solely a function of $r$ and keeping $g$ constant,
we cut off $f_{\exp}$ smoothly starting at some point past $r_c$ to
give a Hamiltonian $H_2$ which agrees with $H_{sing}+g$ outside a
ball. As long as $\frac{\partial f_{exp}}{\partial r}<2\mu r$, there
are no new critical points.

Note that $f_{exp}(r_c) = -\frac{\mu}{m}$. Keeping $\frac{\partial
f_{exp}}{\partial r}$ near $\mu r_{c}$ (which, using e.g. $A
=\frac{e\mu}{m}$, is $\frac{1}{\sqrt{m}}$), we can bring $f_{exp}$
to zero in a radial distance of a constant times
$\frac{1}{\sqrt{m}}$; i.e. for $m$ large we can make $H_2$ agree
with $H_{sing}+g$ outside an arbitrarily small ball.

For $A =\frac{e\mu}{m}$, $g=2\frac{\mu}{m}$. Then keeping $g$ solely
a function of $r$, we cut off $g(r,\theta)$ smoothly starting at
some point past the point where $H_{2}=H_{sing}+g$ to give
Hamiltonian $H_{3}$. As long as $\frac{\partial g}{\partial r}>-2\mu
r$, there are no new critical points. We can make it in such a way
that $H_3$ agrees with $H_{sing}$ outside a small ball.

\item[($3$)] Recall that $H_{3}=H_{sing}+f_{exp}+g$ near the origin and $g(r,\theta)=2\frac{\mu}{m}>0$.
Note that $g(r,\theta)$ is small for large $m$. Now keeping $g$
constant we modify $H_{sing}+f_{exp}+g$ near the origin to give us
$H_4$ which is $Br^2-C$ near the origin (for $B>0$), which
corresponds to the Hamiltonian flow rotating at a constant angular
rate.  Since $\frac{\partial H_3}{\partial
r}=\frac{\partial(H_{sing}+f_{exp})}{\partial r}>0$ for $r < r_c$,
we can patch together $Br^2 - C$ near the origin with $H_2$ outside
a small ball of radius less than $r_c$ in a radially symmetric
manner to get $H_4$ such that $\frac{\partial H_4}{\partial r}>0$
for $r < r_c$ (we do this by choosing $C$ sufficiently large). Note
that $H_{4}$ has a critical point at the origin.

\item[($4$)] Finally, to ensure no fixed points of the time-$1$ flow
of the Hamiltonian vector field of $H$, we let $H$ be $H_4$
multiplied by a radially symmetric function which is $\epsilon$ for
$r < R$ (for $\epsilon$ sufficiently small that the only fixed
points of the time-$1$ flow inside radius $R$ are the critical
points and for $R$ large enough that $H_4$ agrees with $H_{sing}$
for $r > R$) and $1$ for $r
> 2R$. This creates no new fixed points in
the region $R < r < 2R$ because $H_4$ and $\frac{\partial
H_4}{\partial r}$ have the same sign there. Now there are no fixed
points of the time-$1$ flow of the Hamiltonian vector field of $H$,
except for the $nk+1$ critical points of $H$ because outside
radius $R$ there are no compact flow lines.
\end{itemize}
Let $p_1,\dots,p_{nk}$ denote the equally spaced saddle points of
$H$ ordered counterclockwise, i.e., $R_{nk}(p_{i})=p_{i+1}$, where
$R_{nk}$ corresponds to the $\frac{2\pi}{nk}$-rotation around
$p$. We note that $H(p_{s})=0$ for $s=1,\dots,nk$. Hence, by
Morse lemma (arguing the same way as in Lemma 3.2
in~\cite{Golovko}) we get that there is a neighborhood $U_{s}$ of
$p_{s}$ such that $H=axy$ on $U_{s}$, where $s=1,\dots,nk$ and
$a>0$. In addition, observe that $H$ is
$\frac{2\pi}{nk}$-symmetric with respect to $\theta$. Therefore,
$U_{s}$'s together with coordinates $(x,y)$ are
$\frac{2\pi}{nk}$-symmetric with respect to $\theta$, i.e.,
$R_{nk}(U_{s})=U_{s+1}$ and coordinates on $U_{s}$ maps to the
coordinate on $U_{s+1}$. Finally, note
that $H=\tilde{B}r^2-\tilde{C}$ on a neighborhood of the center of
$D(r_{sing})$, which we call $U$, where $\tilde{C}>0$ and
$\tilde{B}$ is a small positive number and hence Hamiltonian flow
rotates at a constant rate near the origin.
\end{proof}

\subsection{Sutured contact solid cylinder}\label{section:construction}
In this section, we construct the sutured contact solid cylinder
that we later will glue to get the sutured contact solid torus with
$2n$ sutures of slope $-\frac{k}{l}$, where $n\in \mathbb N$,
$(k,l)=1$ and $k>l>0$.

Let $\gamma_{p,p_{s}}$ be an embedded curve in $\mathbb R^2$ which
starts at $p$ and ends at $p_{s}$ for $s=1,\dots,nk$. For the time
being, we can think about $\gamma_{p,p_{s}}$ as about the segment
connecting $p$ and $p_{s}$.

We start with the following lemma:

\begin{lemma}\label{le1_formanalcon}
There exists a 1-form $\beta$ on $\mathbb R^2$ satisfying the
following:
\begin{itemize}
\item[($1$)] $d\beta > 0$;
\item[($2$)] its singular foliation given by $ker\,\beta$ has isolated
singularities and no closed orbits;
\item[($3$)] $\beta=\frac{\varepsilon_{c}}{2}r^2d\theta$ on $U$ with respect to
the polar coordinates whose origin is at the center of
$D(r_{sing})$; $\beta=\frac{\varepsilon_{sym}}{2}(xdy-ydx)$ on $U_{s}$
with respect to the coordinates from Lemma~\ref{coordinates}, where
$s\in \{1,\dots,nk\}$; $\beta=\frac{1}{2}r^2d\theta$ on $\mathbb
R^2\setminus D(r_{sing})$ with respect to the polar coordinates
whose origin is at the center of $D(r_{sing})$; here
$0<\varepsilon_{c}\ll \varepsilon_{sym}\ll 1$;
\item[($4$)] the set of hyperbolic points of the singular foliation of
$\beta$ is given by $\{ q_{s} \}_{s=1}^{nk}$ such that $q_{s}$
lies on $\gamma_{p,p_{s}}$ outside of $U_{s}$ and $U$;
\item[($5$)] $\beta$ is $\frac{2\pi}{nk}$-symmetric, i.e.,
$R_{nk}^{\ast}(\beta)=\beta$.
\end{itemize}
\end{lemma}
\begin{proof}
Consider a singular foliation $\mathcal F$ on $\mathbb R^2$ which
satisfies the following:
\begin{itemize}
\item[($1$)] $\mathcal F$ is Morse-Smale and has no closed orbits.

\item[($2$)] The singular set of $\mathcal F$ consists of elliptic points and hyperbolic points.
The elliptic points are the equally spaced saddle points of $H$ and
the center of $D(r_{sing})$. The set of hyperbolic points of the
singular foliation of $\beta$ is given by $\{ q_{s} \}_{s=1}^{nk}$
such that $q_{s}$ lies on $\gamma_{p,p_{s}}$ outside of $U_{s}$ and
$U$.

\item[($3$)] $\mathcal F$ is oriented and for one choice of orientation the flow is transverse to
and exits from $\partial D(r_{sing})$.

\item[($4$)] $\mathcal F$ is $\frac{2\pi}{nk}$-symmetric with respect to $\theta$.

\end{itemize}

Next, we modify $\mathcal F$ near each of the singular points so
that $\mathcal F$ is given by $\beta_{0} = \frac{1}{2}(xdy-ydx)$ on
$U_{s}$ with respect to the coordinates from
Lemma~\ref{coordinates} and $\beta_0 = 2xdy+ydx$ near a hyperbolic
point. On $\mathbb R^2\setminus D(r_{sing})$,
$\beta_{0}=\frac{1}{2}r^2d\theta$ with respect to the polar
coordinates whose origin is at the center of $D(r_{sing})$. In
addition, on $U$, $\beta_{0}=\frac{1}{2}r^2d\theta$ with respect to
the polar coordinates whose origin is at the center of
$D(r_{sing})$. From Lemma~\ref{coordinates} it follows that we can
do it in such a way that the modification of $\mathcal F$ is still
$\frac{2\pi}{nk}$-symmetric. Finally, we get $\mathcal F$ given by
$\beta_{0}$, which satisfies $d\beta_0
> 0$ near the singular points and on $\mathbb R^2\setminus D(r_{sing})$.
Now let $\beta = g\beta_0$, where $g$ is a positive function with
$dg(X) \gg 0$ outside of $U\cup(\cup^{nk}_{s=1}U_{s})\cup(\mathbb
R^2\setminus D(r_{sing}))$,
$g|_{\cup^{nk}_{s=1}U_{s}}=\varepsilon_{sym}$,
$g|_{U}=\varepsilon_{c}$, $g|_{\mathbb R^2\setminus D(r_{sing})}=1$
and $X$ is an oriented vector field for $\mathcal F$ (nonzero away
from the singular points). Here $0<\varepsilon_{c}\ll
\varepsilon_{sym}\ll 1$. Since $d\beta = dg\wedge \beta_{0} +
g\wedge d\beta_{0}$, $dg (X) \gg 0$ guarantees that $d\beta
> 0$.
\end{proof}

\begin{figure}[t]
\includegraphics[width=400pt]{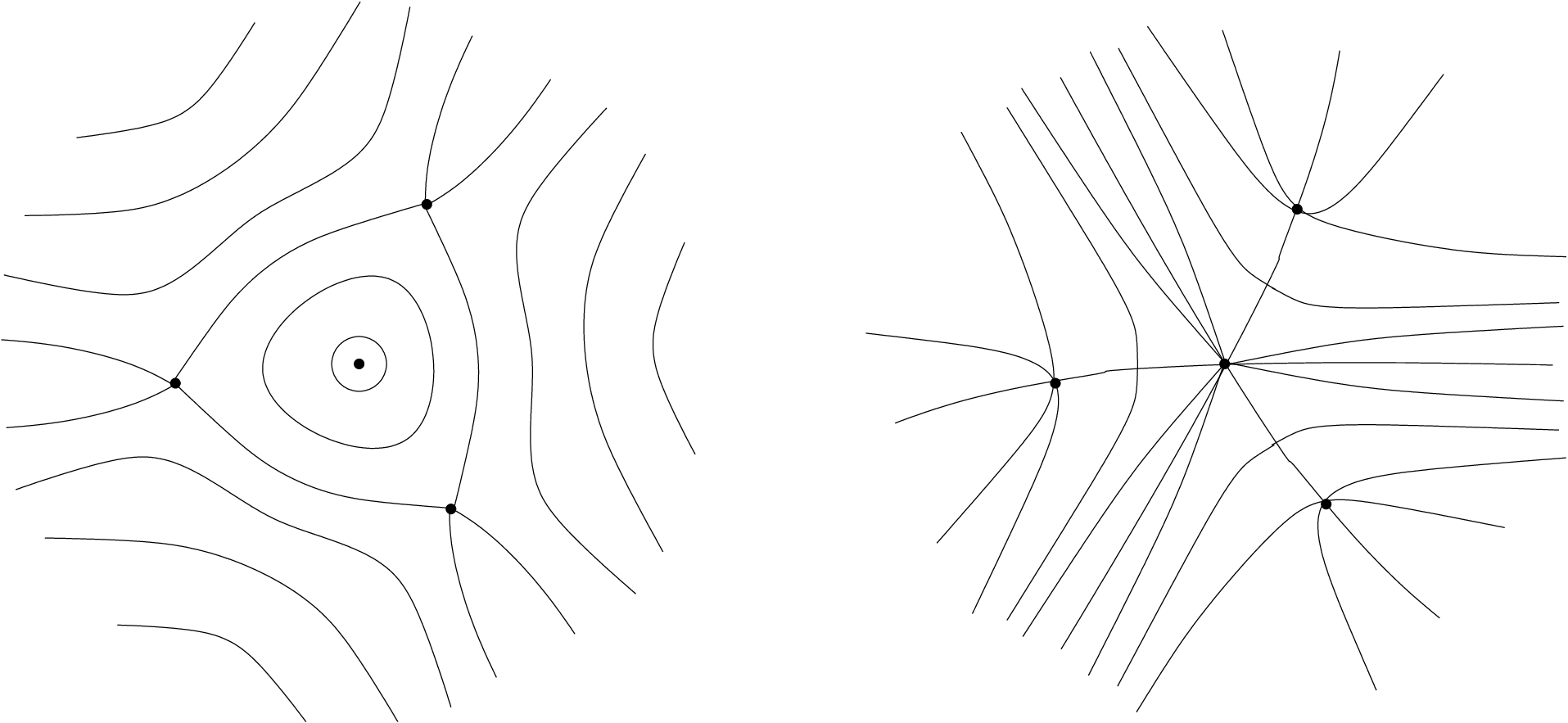}
\caption{The level sets of $H$ (left) and the characteristic
foliation of $\beta$ (right) in the case $n=1$, $k=3$}
\label{charfol2}
\end{figure}

For the comparison of the level sets of $H$ with the singular
foliation of $\beta$ in the case $n=1$, $k=3$ we refer to
Figure~\ref{charfol2}.

\begin{lemma}\label{coordinatesdifference}
Let $\beta$ be a $1$-form from Lemma~\ref{le1_formanalcon}. The
Hamiltonian vector field $X_{H}$ of $H$ with respect to the area
form $d\beta$ satisfies $\beta(X_{H})=H$ on $(\cup_{s=1}^{nk}
U_{s})\cup ( \mathbb R^2\setminus D(r_{sing}))$. In addition, the
Hamiltonian vector field $X_{H}$ of $H$ with respect to the area
form $d\beta$ satisfies $\beta(X_{H})-H=\tilde{C}$ on $U$.
\end{lemma}
\begin{proof}
First, Lemmas~\ref{coordinates} and \ref{le1_formanalcon} imply that
$\beta=\frac{\varepsilon_{c}}{2}r^2d\theta$,
$H=\tilde{B}r^2-\tilde{C}$ on $U$ and $\varepsilon_{c}$ is a small
positive number. Now we show that
$X_{H}=\frac{2\tilde{B}}{\varepsilon_{c}}\frac{\partial}{\partial
\theta}$ is a solution of $\beta(X_{H})-H=\tilde{C}$ on $U$. We
calculate
\begin{align*}
i_{X_{H}}(d\beta)&=\left(\frac{2\tilde{B}}{\varepsilon_{c}}\frac{\partial}{\partial
\theta}\right)\lrcorner(\varepsilon_{c}rdr\wedge
d\theta)=-2\tilde{B}rdr=-dH,
\end{align*}
and
\begin{align*}
\beta(X_{H})-H=\left (\frac{\varepsilon_{c}}{2}r^2d\theta \right
)\left(\frac{2\tilde{B}}{\varepsilon_{c}}\frac{\partial}{\partial
\theta}\right)-\tilde{B}r^2+\tilde{C}=\tilde{C}.
\end{align*}

Next, we work on $U_{s}$, where $s=1,\dots,nk$. From
Lemmas~\ref{coordinates} and \ref{le1_formanalcon} it follows that
$\beta=\frac{\varepsilon_{sym}}{2}(xdy-ydx)$ and $H=axy$ on $U_{s}$.
Let $X_{H}$ be a Hamiltonian vector field defined by $i_{X_{H}}
d\beta=-dH$.

We show that
\begin{align*}X_{H}=-\frac{ax}{\varepsilon_{sym}}
\frac{\partial}{\partial x}+\frac{ay}{\varepsilon_{sym}}
\frac{\partial}{\partial y}
\end{align*}
is a solution of the equation
\begin{align}\label{diff_eq1}
\beta(X_{H})=H
\end{align}
on $U_{s}$. We calculate
\begin{align*}
i_{X_{H}}( d\beta)=\left(-\frac{ax}{\varepsilon_{sym}}
\frac{\partial}{\partial x}+\frac{ay}{\varepsilon_{sym}}
\frac{\partial}{\partial y}\right)\lrcorner
(\varepsilon_{sym}dx\wedge dy)=-axdy-aydx=-dH
\end{align*}
and
\begin{align*}
\beta(X_{H})=\frac{\varepsilon_{sym}}{2}(xdy-ydx)\left
(-\frac{ax}{\varepsilon_{sym}} \frac{\partial}{\partial
x}+\frac{ay}{\varepsilon_{sym}} \frac{\partial}{\partial y}\right
)=axy=H.
\end{align*}

Finally, Lemmas~\ref{coordinates} and \ref{le1_formanalcon} say that
$\beta=\frac{1}{2}r^2d\theta$ and $H=\mu r^2\cos(nk\theta)$ on
$\mathbb R^2\setminus D(r_{sing})$. As in the previous case, we show
that
\begin{align*}
X_{H}=nk\mu r\sin(nk\theta)\frac{\partial}{\partial r}+2\mu
\cos(nk\theta)\frac{\partial}{\partial \theta}
\end{align*}
is a solution of Equation~(\ref{diff_eq1}) on $\mathbb R^2\setminus
D(r_{sing})$.

We calculate
\begin{align*}
i_{X_{H}}(d\beta)&=(nk\mu r \sin(nk\theta)\partial_{r}+2\mu
\cos(nk\theta)\partial_{\theta})\lrcorner(rdr\wedge d\theta)\\
&=-2\mu r \cos(nk\theta)dr+nk\mu r^2\sin(nk\theta)d\theta=-dH,
\end{align*}
and
\begin{align*}
\beta(X_{H})&=\left (\frac{1}{2}r^2d\theta \right ) \left (nk\mu r
\sin(nk\theta)\frac{\partial}{\partial r}+2\mu
\cos(nk\theta)\frac{\partial}{\partial \theta}\right )\\ &=\mu r^2
\cos(nk\theta)=H.
\end{align*}

\end{proof}

Let $X_{H}$ be the Hamiltonian vector field of $H$ with respect to
$d\beta$ and $\varphi^{s}_{X_{H}}$ be the time-$s$ flow of $X_{H}$.
Now we introduce the following notations:
\begin{align*}
&S:=\{ x\in  \mathbb R^2\setminus D(r_{sing})\ \vert \
\varphi_{X_{H}}^{s}(x)\in \mathbb R^2\setminus D(r_{sing})\
\forall s\in[0,1] \},\\
&V:=\{x\in U\ \vert \ \varphi_{X_{H}}^{s}(x)\in U\ \forall
s\in[-1,1] \},\quad
\mbox{and}\\
&V_{i}:=\{ x\in  U_{i}\ \vert \ \varphi_{X_{H}}^{s}(x)\in U_{i}\
\forall s\in[-1,1] \}.
\end{align*}
For simplicity, let us denote
$\varphi_{X_{H}}:=\varphi_{X_{H}}^{1}$.

\begin{remark}\label{extension}
Using the form of $X_{H}$ on $U_{i}$, where $i=1,\dots,nk$, we may
assume that the curves $\gamma_{p,p_{i}}$'s in
Lemma~\ref{le1_formanalcon} satisfy the following list of properties:
\begin{itemize}
\item[($1$)] $\gamma_{p,p_{i}}$ is an embedded curve which starts at $p$ and ends at $p_{i}$;
\item[($2$)] $\gamma_{p,p_{i}}$ is a part of one of the curves of the singular foliation given
by $ker\,\beta$;
\item[($3$)] $\gamma_{p,p_{i}}$ coincides with one of the level sets of $H$ on $V_{i}$ and near $p_{i}$ can be presented as $W^{s}(\varphi_{X_{H}},p_{i})=\{x\
\vert\ (\varphi_{X_{H}})^{n}(x)\rightarrow p\ \mbox{as}\
n\rightarrow \infty \}$.
\end{itemize}
\end{remark}

Recall that the following claim was proven in~\cite{Golovko}:
\begin{claim}[\cite{Golovko}]\label{difference}
If $(M,\omega)$ is an exact symplectic manifold, i.e.,
$\omega=d\beta$, then the flow $\varphi_{X_{H}}^{t}$ of a
Hamiltonian vector field $X_{H}$ consists of exact symplectic maps,
i.e.,
\begin{align*}
(\varphi_{X_{H}}^{t})^{\ast}\beta - \beta = df_{t},
\end{align*}
where
\begin{align*}
f_t=\int\limits^{t}_{0}(-H+\beta(X_{H}))\circ \varphi_{X_{H}}^{s}ds.
\end{align*}
\end{claim}

\begin{remark}
Observe that from Lemma~\ref{coordinatesdifference} and
Claim~\ref{difference} it follows that
$\varphi_{X_{H}}^{\ast}(\beta)-\beta=dh$, where $h:=f_{1}=0$ on
$S\cup(\cup_{i=1}^{nk} V_{i})$ and $h=\tilde{C}>0$ on $V$. Hence,
we get $\varphi_{X_{H}}^{\ast}(\beta)=\beta$ on $S\cup
V\cup(\cup_{i=1}^{nk} V_{i})$.
\end{remark}

Now we define $\varphi_{-\frac{k}{l}}:=R_{-\frac{k}{l}}\circ
\varphi_{X_{H}}$, where $R_{-\frac{k}{l}}:\mathbb R^2 \rightarrow
\mathbb R^2$ is a $-\frac{2\pi l}{k}$-rotation around $p$.

\begin{remark}
Since $R_{nk}^{\ast}(\beta)=\beta$, we get
$R_{-\frac{k}{l}}^{\ast}(\beta)=\beta$ and hence
\begin{align*}
\varphi_{-\frac{k}{l}}^{\ast}(\beta)=(R_{-\frac{k}{l}}\circ
\varphi_{X_{H}})^{\ast}(\beta)=\varphi_{X_{H}}^{\ast}(R_{-\frac{k}{l}}^{\ast}(\beta))=\varphi_{X_{H}}^{\ast}(\beta).
\end{align*}
\end{remark}

Fix $R_{\ast}\gg r_{sing}$ such that there is an annular
neighborhood $V_{R_{\ast}}$ of $\partial D(R_{\ast})$ in $\mathbb
R^2$ with $V_{R_{\ast}}\subset S$. Consider $D(R_{\ast})$ with
$\beta_{0}:=\beta|_{D(R_{\ast})}$ and
$\beta_{1}:=\varphi_{X_{H}}^{\ast}(\beta)|_{D(R_{\ast})}(=\varphi_{-\frac{k}{l}}^{\ast}(\beta)|_{D(R_{\ast})})$.
Note that
\begin{align}\label{arform}
d\beta_{1}=d(\varphi_{X_{H}}^{\ast}(\beta)|_{D(R_{\ast})})=\varphi_{X_{H}}^{\ast}(d\beta)|_{D(R_{\ast})}=
(d\beta)|_{D(R_{\ast})}=d\beta_{0}>0.
\end{align}
In addition, from the definitions of $V(R_{\ast})$ and $D(R_{\ast})$
it follows that
\begin{align}\label{boundaryeql}
\beta_{0}=\beta_{1}\quad \mbox{on}\quad V_{R_{\ast}}\cap
D(R_{\ast}).
\end{align}

Now we recall Lemma~3.10 from~\cite{Golovko}, which provides the
construction of the contact $1$-form on $[-1,1]\times D^2$.
\begin{lemma}[\cite{Golovko}]\label{befglcontformconstr}
Let $\beta_0$ and $\beta_1$ be $1$-forms on $D^2$ such that
$\beta_0=\beta_1$ in a neighborhood of $\partial D^2$ and
$d\beta_0=d\beta_1=\omega>0.$ Then there is a contact $1$-form
$\alpha$ and a Reeb vector field $R_{\alpha}$ on $[-1,1] \times D^2$
with coordinates $(t,x)$, where t is a coordinate on $[-1,1]$ and
$x$ is a coordinate on $D^2$, such that:
\begin{itemize}
\item[($1$)] $\alpha=dt+\varepsilon\beta_{0}$ in a neighborhood of $\{-1\}\times
D^2$;
\item[($2$)] $\alpha=dt+\varepsilon\beta_{1}$ in a neighborhood of $\{1\}\times
D^2$;
\item[($3$)] $R_{\alpha}$ is collinear to $\frac{\partial}{\partial
t}$ on $[-1,1]\times D^2$;
\item[($4$)] $R_{\alpha}=\frac{\partial}{\partial t}$ in a
neighborhood of $[-1,1]\times \partial D^2$.
\end{itemize}
Here $\varepsilon$ is a small positive number.
\end{lemma}
In addition, recall that
\begin{align}\label{befglcontform}
\alpha=(1+\varepsilon
\chi_1(t)h)dt+\varepsilon((1-\chi_{0}(t))\beta_{0}+\chi_{0}(t)\beta_{1}),
\end{align}
where $h \in C^{\infty}(D^2)$ such that $\beta_1 - \beta_0 = dh$;
$\chi_{0}: [-1,1]\rightarrow [0,1]$ is a smooth map for which
$\chi_{0}(t)=0$ for $-1\leq t \leq-1+\varepsilon_{\chi_{0}}$,
$\chi_{0}(t)=1$ for $1-\varepsilon_{\chi_{0}}\leq t \leq 1$,
$\chi_{0}'(t)\geq 0$ for $t\in [-1,1]$ and  $\varepsilon_{\chi_{0}}$
is a small positive number; $\chi_{1}(t):=\chi_{0}'(t)$;
$\varepsilon$ is a sufficiently small positive number.

\begin{remark}\label{approxform}
Note that $d\alpha=\varepsilon \omega$, where $\alpha$ is a $1$-form
given by Formula~\ref{befglcontform} and
$\omega=d\beta_{0}=d\beta_{1}>0$ on $D^2$.
\end{remark}

Observe that from Formulas~\ref{arform} and \ref{boundaryeql} it
follows that $\beta_{0}$ and $\beta_{1}$ described above satisfy the
conditions of Lemma~\ref{befglcontformconstr}. We now take
$[-1,1]\times D(R_{\ast})$ equipped with the contact $1$-form
$\alpha$ given by Formula~\ref{befglcontform}. For simplicity, let
us denote $\beta_{-}:=\varepsilon\beta_{0}$ and
$\beta_{+}:=\varepsilon\beta_{1}$, where $\varepsilon$ is a constant
from Lemma~\ref{befglcontformconstr} which makes $\alpha$ contact.

\begin{figure}[t]
\begin{center}
\labellist
\pinlabel $D(r_{sing})$ at 110 105
\pinlabel $a^{+}_{0}$ at 110 190
\pinlabel $a^{+}_{1}$ at 25 55
\pinlabel $a^{+}_{2}$ at 195 55
\pinlabel $b^{-}_{0}$ at 25 155
\pinlabel $b^{-}_{1}$ at 110 15
\pinlabel $b^{-}_{2}$ at 194 155
\pinlabel $c^{+}_{0}$ at -7 110
\pinlabel $c^{+}_{1}$ at 167 0
\pinlabel $c^{+}_{2}$ at 165 208
\pinlabel $c^{-}_{0}$ at 60 209
\pinlabel $c^{-}_{1}$ at 50 0
\pinlabel $c^{-}_{2}$ at 227 110
\endlabellist
\includegraphics[width=215pt]{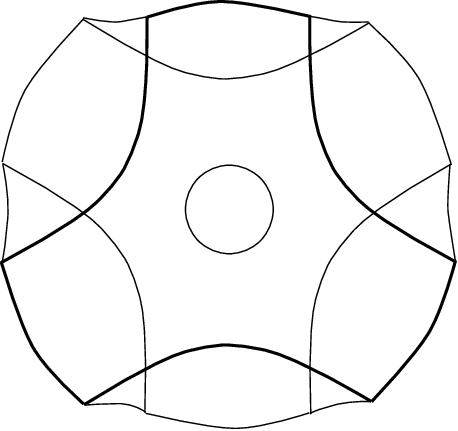}
\caption{Construction of $P_{+}$, $P_{-}$ and $D$ in the case $n=1$,
$k=3$.}\label{constofd}
\end{center}
\end{figure}
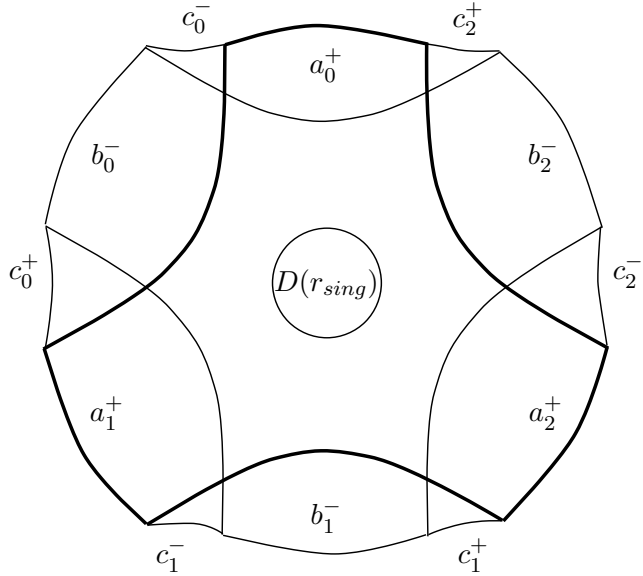

\subsection{Gluing}

We now construct $P_{+}$, $P_{-}$ and $D$ in the way described
in~\cite{Golovko}. Recall that
\begin{align*}
P_{+},P_{-}, D\subset D(R_{\ast}) \subset \mathbb R^2
\end{align*}
are surfaces with boundary which satisfy the following properties:
\begin{itemize}
\item[($1$)] $P_{\pm}\subset D$;
\item[($2$)] $(\partial P_{\pm})_{\partial} \subset \partial D$ and $(\partial P_{\pm})_{int}\subset int(D)$;
\item[($3$)] $\varphi_{X_{H}}$ maps $P_{+}$ to $P_{-}$ in such a way that $\varphi_{X_{H}}((\partial P_{+})_{int})=(\partial P_{-})_{\partial}$ and
$\varphi_{X_{H}}((\partial P_{+})_{\partial})=(\partial
P_{-})_{int}$;
\item[($4$)] $(\partial P_{-})_{\partial} \cap (\partial P_{+})_{\partial} =
\emptyset$.
\end{itemize}
Note that 
\begin{itemize}
\item $\partial P_{+}=(\cup^{nk-1}_{s=0} a^{+}_{s})\cup
(\cup^{nk-1}_{s=0} b^{+}_{s})$,

\item $\partial
P_{-}=(\cup^{nk-1}_{s=0} a^{-}_{s})\cup (\cup^{nk-1}_{s=0}
b^{-}_{s})$,

\item $\partial D=(\cup^{nk-1}_{s=0} a^{+}_{s})\cup
(\cup^{nk-1}_{s=0} b^{-}_{s})\cup (\cup^{nk-1}_{s=0}
c^{+}_{s})\cup (\cup^{nk-1}_{s=0} c^{-}_{s})$.
\end{itemize}

See Figure~\ref{constofd} for the schematic visualization of $P_{+}$
(bounded by the bold line), $P_{-}$ and $D$. For more details of
this construction we refer to~\cite{Golovko}.

\begin{remark}\label{relntos}
Note that $a^{\pm}_{i}$'s, $b^{\pm}_{i}$'s and $c^{\pm}_{i}$'s are
constructed in such a way that
\begin{align*}
a^{\pm}_{i}, b^{\pm}_{i}, c^{\pm}_{i}\subset D(R_{\ast})\cap S
\end{align*}
for $i=0,\dots,nk-1$. Hence, we see that $\partial P_{+},
\partial P_{-},
\partial D\subset D(R_{\ast})\cap S$. In addition,
$R_{nk}(a^{\pm}_{i})=a^{\pm}_{i+1}$ and
$R_{nk}(b^{\pm}_{i})=b^{\pm}_{i+1}$, where $i,i+1$ are considered
modulo $nk$.
\end{remark}

We take $[-1,1]\times D$ with a contact form
$\alpha:=\alpha|_{[-1,1]\times D}$. Let $\Gamma=\{0\}\times \partial
D$ in $[-1,1]\times D$ and $U(\Gamma):=[0,1]\times [-1,1]\times
\Gamma$ be a neighborhood of $\Gamma$ with coordinates $(\tau,t)\in
[0,1]\times [-1,1]$, where $t$ is a usual $t$-coordinate on
$[-1,1]\times D$. From the definition of $S$ and
Remark~\ref{relntos} it follows that we may assume that
$U(\Gamma)\subset [-1,1]\times (S\cap D)$.

\begin{lemma}\label{le1_suturbefglug}
$([-1,1]\times D, \Gamma, U(\Gamma), \xi)$ is a sutured contact
manifold and $\alpha$ is an adapted contact form.
\end{lemma}
\begin{proof}
First note that $\alpha|_{R{-}}=\beta_{-}$ and
$\alpha|_{R_{+}}=\beta_{+}$. Let us check that $(R_{-},\beta_{-})$
and $(R_{+},\beta_{+})$ are Liouville manifolds. From the
construction of $\beta_{\pm}$ it follows that
$d(\beta_{-})=d(\beta_{+})>0$. Since $\beta_{-}=\beta_{+}$ on $D\cap
S$ and by Formula~\ref{befglcontform}, $\alpha=dt+\beta_{-}$ on
$U(\Gamma)$. Recall that
$\beta_{-}=\beta_{+}=\frac{\varepsilon}{2}r^2d\theta$ on $D\cap S$.
Hence, $\alpha|_{U(\Gamma)}=dt+\frac{\varepsilon}{2} r^2d\theta$.
The calculation
\begin{align*}
i_{Y_{\pm}|_{R_{\pm}\cap
U(\Gamma)}}(d\beta_{\pm})=\left(\frac{1}{2}r\partial_{r}\right)\lrcorner(\varepsilon
rdr\wedge d\theta)=\frac{\varepsilon}{2} r^2d\theta=\beta_{\pm}
\end{align*}
implies that the Liouville vector fields $Y_{\pm}|_{R_{\pm}\cap
U(\Gamma)}$ are equal to $\frac{1}{2}r\partial_{r}$. From the
construction of $D$ it follows that $Y_{\pm}$ is positively
transverse to $\partial R_{\pm}$. Therefore,
$(R_{-},\varepsilon\beta_{0})$ and $(R_{+},\varepsilon\beta_{1})$
are Liouville manifolds. As we already mentioned,
$\alpha=dt+\beta_{-}$ on $U(\Gamma)$. Finally, if we take $\tau$
such that $\partial_{\tau}=\frac{1}{2}r\partial_{r}$, then
$([-1,1]\times D, \Gamma, U(\Gamma), \xi)$ becomes a sutured contact
manifold with an adapted contact form $\alpha$.
\end{proof}

Then we use $\varphi_{-\frac{k}{l}}$ for the gluing construction.
Note that $\varphi_{X_{H}}$ maps $a^{+}_{s}$ to $a^{-}_{s}$ and
$b^{+}_{s}$ to $b^{-}_{s}$. Hence, using Remark~\ref{relntos}, we
see that $\varphi_{-\frac{k}{l}}$ maps $a^{+}_{s}$ to $a^{-}_{s-nl}$
and $b^{+}_{s}$ to $b^{-}_{s-nl}$. Then we follow the gluing
procedure briefly described in Section~\ref{section:glsutcontman}
and completely written in~\cite{ColinGhigginiHondaHutchings}.
Finally, we get a sutured contact solid torus $(S^1\times D^2,
\tilde{\Gamma}, U(\tilde{\Gamma}))$ with a contact form
$\tilde{\alpha}_{\delta}$, where $\tilde{\Gamma}$ is a set of $2n$
parallel closed curves of slope $-\frac{k}{l}$, where $n\in \mathbb
N$, $(k,l)=1$, $k>l>0$ and $\delta$ is the rotation angle of the map
$\varphi_{X_{H}}$ near $p$.

\begin{remark}
We have constructed $(S^1\times D^2,\tilde{\Gamma},
U(\tilde{\Gamma}))$ using the gluing construction for sutured
manifolds. However, since there is a close connection between
sutured contact manifolds and contact manifolds with convex
boundary, we observe that the gluing construction we used for the
sutured contact solid cylinder corresponds to the gluing
construction for the contact $3$-ball with convex boundary and one
dividing curve on the boundary. The corresponding gluing
construction for the contact $3$-ball with convex boundary
corresponds (is inverse) to the convex decomposition of the contact
solid torus $S^1\times D^2$ with convex boundary with respect to the
convex meridional disk $\{ pt \}\times D^2$ with $\partial$-parallel
dividing curves. Hence, the constructed sutured contact solid tori
are universally tight sutured contact manifolds by the gluing/classification result from Section 2 in~\cite{Honda} (more precisely, Corollary~2.3, Theorem~2.5 and Corollary~2.6).
\end{remark}

\subsection{Reeb orbits} Note that $\varphi_{-\frac{k}{l}}|_{P_{+}}$ has $n$
orbits of period $k$ obtained from the equally spaced saddle
points of $H$. Lemma~\ref{befglcontformconstr} and the gluing
procedure briefly described in Section~\ref{section:glsutcontman}
imply that these orbits correspond to the Reeb orbits, which we call
$\gamma_1,\dots,\gamma_n$ such that
\begin{align*}
[\gamma_{s}]=[\gamma_{t}]=k\in H_{1}(S^1\times D^2;\mathbb Z)
\end{align*}
for $s,t=1,\dots,n$. In addition, $\varphi_{-\frac{k}{l}}|_{P_{+}}$
has a periodic point of period $1$, which is $p$. It corresponds to
the Reeb orbit, which we call $\gamma$, such that $[\gamma]=1\in
H_{1}(S^1\times D^2;\mathbb Z)$.

\begin{lemma}\label{orbits}
$\int_{\gamma_{s}}\tilde{\alpha}_{\delta}=\int_{\gamma_{t}}\tilde{\alpha}_{\delta}$
and
$k\int_{\gamma}\tilde{\alpha}_{\delta}>\int_{\gamma_{s}}\tilde{\alpha}_{\delta}$,
where $s,t=1,\dots,n$.
\end{lemma}
\begin{proof}
Let
\begin{align*}
M^{(0)}=(([-1,1]\times D)\cup (R_{+}(\Gamma)\times[1;\infty))\cup
(R_{+}(\Gamma)\times (-\infty;-1]))
\end{align*}
and
\begin{align*}
\tilde{M}=M^{(0)} \setminus ((P_{+}\times (N, \infty)\cup
(P_{-}\times (-\infty,-N)).
\end{align*}
In addition, let $\alpha_{\tilde{M}}$ denote the contact form on
$\tilde{M}$ and let $\xi_{\tilde{M}}$ denote the contact structure
defined by $\alpha_{\tilde{M}}$.

Consider $[-1,1]\times D\subset \tilde{M}$. From the construction of
$\alpha$ it follows that $\beta_{+}=\beta_{-}$ on $V_{s}$ and
$\alpha|_{[-1,1]\times V_{s}}=dt+\beta_{-}$ for $s=1,\dots,nk$.
Hence, since the contact structure on $[1,\infty)\times P_{+}$ is
given by $dt+\beta_{+}$ and the contact structure on $(-\infty,
-1]\times P_{-}$  is given by $dt+\beta_{-}$,
$\alpha_{\tilde{M}}|_{[-N,N]\times V_{s}}=dt+\beta_{-}$ on
$[-N,N]\times V_{s}\subset \tilde{M}$ for $s=1,\dots,nk$.
Therefore, we get
\begin{align}\label{e10_anconstlength}
\int\limits_{[-N,N]\times \{p_{s}\}} \alpha_{\tilde{M}}=2N
\end{align}
for $s=1,\dots,nk$. From the gluing construction and
Equation~(\ref{e10_anconstlength}) it follows that
\begin{align*}
\int\limits_{\gamma_{s}} \tilde{\alpha}_{\delta}=2Nk
\end{align*}
for $s=1,\dots,n$. Note that $\int_{\gamma_{s}}
\tilde{\alpha}_{\delta}$ does not depend on $s$. Hence,
$\int_{\gamma_{s}}\tilde{\alpha}_{\delta}=\int_{\gamma_{t}}\tilde{\alpha}_{\delta}$
for $s,t=1,\dots,n$.

Now from the fact that $\alpha=(1+\varepsilon
\chi_1(t)h)dt+\beta_{-}$ on $[-1,1]\times V$, where $h>0$ and
$\chi_1(t)>0$, we get that
\begin{align*}
R_{\alpha}=\frac{1}{1+\varepsilon \chi_1(t)h}
\frac{\partial}{\partial t}
\end{align*}
on $[-1,1]\times V$. Hence, from the gluing construction we obtain
$k\int_{\gamma} \tilde{\alpha}_{\delta}>2Nk$. Thus,
\begin{align*}
\int\limits_{\gamma_{s}}\tilde{\alpha}_{\delta}=\int\limits_{\gamma_{t}}\tilde{\alpha}_{\delta}\quad
\mbox{and}\quad
k\int\limits_{\gamma}\tilde{\alpha}_{\delta}>\int\limits_{\gamma_{s}}\tilde{\alpha}_{\delta},
\end{align*}
where $s,t=1,\dots,n$.
\end{proof}

\begin{lemma}\label{trivialization}
All closed orbits of $R_{\tilde{\alpha}_{\delta}}$ are
nondegenerate. Moreover, $\gamma$ is an elliptic orbit and
$\gamma_{i}$ is a hyperbolic orbit such that $\gamma^{t}$ and
$\gamma^{s}_{i}$ are good orbits for $i=1,\dots n$; $s,t\in \mathbb
N$. There exists a symplectic trivialization $\tau$ of $\xi$ along
$\gamma$ and $\gamma_{i}$'s constructed in the consistent way as
described in Section~\ref{section:sutcylconthom}, and $N_{\delta}\in
\mathbb N$ such that

\begin{align*}
&\mu_{\tau}(\gamma^{s}_{i})=
-2ls\quad \mbox{and}\\
&\mu_{\tau}(\gamma^t)=-2ml+1,
\end{align*}
where $(m-1)k< t\leq mk$ and $i=1,\dots,n$, $t\leq N_{\delta}$,
$s\leq \frac{N_{\delta}}{k}$.
\end{lemma}

\begin{proof}
For simplicity, assume that $l=1$. The general calculation can be done in the analogous way.

Fix $i=1,\dots,n$. We first observe that $H|_{V_{i}}=axy$, where
$a>0$ and hence
\begin{align*}\varphi_{X_{H}}|_{V_{i}}=
\left ( \begin{array}{ll} \lambda & 0\\
0 & \lambda^{-1}
\end{array}
\right ),
\end{align*}
where $\lambda=e^{a}\neq 1$. Let the symplectic trivialization of
$\xi_{\tilde{M}}$ along $[-N,N]\times \{p_{i}\}$ be given by the
framing
$(\lambda^{\frac{-N-t}{2N}}\partial_{x},\lambda^{\frac{t+N}{2N}}\partial_{y})$,
where $i=1,\dots,nk$ and $(x,y)$ are coordinates on $V_{i}$ which
coincide with the coordinates on $U_{i}$ from
Lemma~\ref{coordinates}. Since Lemma~\ref{coordinates} implies that
$R_{nk}$ maps coordinates on $V_{i}$ to the coordinate on
$V_{i+1}$, where $i$, $i+1$ are considered modulo $nk$, we
conclude that the symplectic trivializations of $\xi_{\tilde{M}}$
along $[-N,N]\times \{p_{i+nm}\}$'s for $m=0,\dots,k-1$ and fixed
$i=1,\dots,n$ give rise to the symplectic trivialization
$\tau_{\gamma_{i}}$ of $\tilde{\xi}$ along $\gamma_{i}$. It is easy
to see that the linearized return map $P_{\gamma_{i}}$ with respect
to this trivialization is given by
\begin{align*}P_{\gamma_{i}}=
\left ( \begin{array}{ll} \lambda^{k} & 0\\
0 & \lambda^{-k}
\end{array}
\right ).
\end{align*}

Since the eigenvalues of $P_{\gamma_{i}}$ are positive real numbers
different from $1$, $\gamma_{i}$ is a positive hyperbolic orbit. In
addition, $P_{\gamma^{s}_{i}}=P^s_{\gamma_{i}}$. Therefore, the
eigenvalues of $P_{\gamma^{s}_{i}}$ are different from $1$. Hence,
$\gamma^{s}_{i}$ is a nondegenerate orbit for $s\in \mathbb N$ and
$i=1,\dots,n$. We now observe that the linearized Reeb flow around
$\gamma_{i}$ (with respect to $\tau_{\gamma_{i}}$) rotates the
eigenspaces of $P_{\gamma_{i}}$ by angle $-2\pi$. Hence, we get
\begin{align}\label{lochyptriv}
\mu_{\tau_{\gamma_{i}}}(\gamma^{s}_{i})= -2s
\end{align}
for $s\in \mathbb N$ and $i=1,\dots,n$.

Now let the symplectic
trivialization of $\xi_{\tilde{M}}$ along $[-N,N]\times \{p\}$ be
given by the framing
\begin{align*}
(\cos(\theta_{\delta,k,N}(t))\partial_{x}+\sin(\theta_{\delta,k,N}(t))\partial_{y},
-\sin(\theta_{\delta,k,N}(t))\partial_{x}+\cos(\theta_{\delta,k,N}(t))\partial_{y}),
\end{align*}
where $\theta_{\delta,k,N}(t)=\frac{\pi(1-\delta k)(t+N)}{Nk}$ and
$t\in [-N,N]$. Observe that $R_{-k}\circ \varphi_{X_{H}}|_{V}$ is a
rotation by angle $2\pi(-\frac{1}{k}+\delta)$, where $R_{-k}$ is a
$-\frac{2\pi}{k}$-rotation about $p$ and $\delta$ is a small
positive irrational number. It is easy to see that with respect to
this framing $P_{\gamma}$ is a rotation by
$2\pi(-\frac{1}{k}+\delta)$. Hence, since $\delta$ is irrational, we
see that $\gamma$ is an elliptic orbit and $\gamma^{t}$ is
nondegenerate for $t\in \mathbb N$. Let
\begin{align*}
N_{\delta}:= \max\{m\in \mathbb N \ \vert\ m\delta<\frac{1}{k} \}.
\end{align*}
Note that we get
\begin{align}\label{loceltriv}
\mu_{\tau_{\gamma}}(\gamma^t)= -2m+1,
\end{align}
where $(m-1)k<t\leq mk$ and $t\leq N_{\delta}$.
Formulas~\ref{lochyptriv} and~\ref{loceltriv} and the fact that
$\delta$ is irrational imply that the parity of
$\mu_{\tau_{\gamma_{i}}}(\gamma^{s}_{i})$ is independent of $s$ for
given $i$ and the parity of $\mu_{\tau_{\gamma}}(\gamma^{t})$ is
independent of $t$. Hence, we conclude that $\gamma^{s}_{i}$'s and
$\gamma^{t}$'s are good Reeb orbits for $i=1,\dots,n$ and $s,t\in
\mathbb N$.

It is not difficult to see that the symplectic trivialization
$\tau_{\gamma^{k}}$ (induced from $\tau_{\gamma}$) can be extended
to $\tau_{\gamma_i}$'s (are consistent in terms of
Section~\ref{section:sutcylconthom}) along the surfaces obtained
from $(\varphi^{(-N-t)/2N}_{X_{H}}(\gamma_{p,p_{i}}))_{i=1}^{nk}$ by
gluing them with $\varphi_{-k}$ and gives rise to the global
symplectic trivialization that we call $\tau$.

\end{proof}

\section{Calculation}

In this section, we calculate the sutured
version of cylindrical contact homology of the sutured contact solid
torus that we have constructed in
Section~\ref{constrodthesoltornaint}.

\begin{remark}\label{cch_indep}
Note that there are no contractible Reeb orbits. Hence, from
Theorem~\ref{welldefconthom}, Remark~\ref{translationforsutures},
and the fact that $\pi_{1}(S^1\times D^2; \mathbb Z)\simeq
H_{1}(S^1\times D^2; \mathbb Z)\simeq \mathbb Z$ it follows that for
all $h\in H_{1}(S^1\times D^2; \mathbb Z)$, $HC^{cyl,
h}_{\ast}(S^1\times D^2,\tilde{\Gamma},\tilde{\alpha}_{\delta}, J)$
is defined, i.e., $\partial^2=0$, and is independent of contact form
$\tilde{\alpha}_{\delta}$ for the given contact structure
$\tilde{\xi}$ and the almost complex structure $J$.
\end{remark}

For simplicity, assume that $l=1$. The calculation for $l>1$ can be made in the completely analogous way.

Lemma~\ref{trivialization} implies that all Reeb orbits are good and
\begin{align}\label{cz1}
&|\gamma_i^{s}|=
-2s-1,\\
&|\gamma^t|= -2m, \nonumber
\end{align}
where $(m-1)< \frac{t}{k}\leq m$ and $i=1,\dots,n$, $s\leq
\frac{N_{\delta}}{k}$, $t\leq N_{\delta}$. Hence, we get
\begin{equation}\label{chaincomplex1}
C_{m}^{h}(\tilde{\alpha_{\delta}}, J)=\left \{
\begin{array}{ll}
\mathbb Q \langle \gamma^{h} \rangle, & \mbox{for}\ h > 0\ \mbox{and} \ m=2\lfloor h(-\frac{1}{k}+\delta) \rfloor; \\
\mathbb Q \langle \gamma_{1}^{h/k},\dots,\gamma_{n}^{h/k}
\rangle, & \mbox{for}\ k \mid h > 0\ \mbox{and}\ m=-\frac{2h}{k}-1;\\
0, & \mbox{otherwise}
\end{array}
\right.
\end{equation}
for $h\leq N_{\delta}$.

Now, since by Lemma~\ref{orbits} $\mathcal A(\gamma^{k})
> \mathcal A(\gamma_{i})$ for $i=1,\dots,n$, we can use
Fact~\ref{energy} and Remark~\ref{translationforsutures} and
conclude that $\partial (\gamma^{s}_{i})=0$ for $i=1,\dots,n$ and
$s>0$. Then, we prove that $\partial (\gamma^{t})=0$ for $k \nmid
t\leq N_{\delta}$. Since $[\gamma_{i}]=k[\gamma]$ in
$H_{1}(S^1\times D^2; \mathbb Z)\cong \mathbb Z$, the cylindrical
contact homology differential at $\gamma^{t}$ counts only cylinders
with negative end at $\gamma^{t}$. Then, similarly to the previous
case, Fact~\ref{energy} and Remark~\ref{translationforsutures} imply
that $\partial(\gamma^{t})=0$ for $k \nmid t\leq N_{\delta}$.

We now consider the case when $k \mid t$ and will show that $\partial(\gamma^{t})\neq 0$ for $k\mid t\leq N_{\delta}$.
Is this situation, by
arguing in the same way as in the case when $k \nmid t$, we get that
$\partial(\gamma^{t})$ counts only cylinders with negative end at
$\gamma_{i}^{t/k}$.

Now we note
that
\begin{align}\label{index}
ind(u)=|\gamma^{t}|-|\gamma_{i}^{t/k}|
\end{align}
for any pseudoholomorphic curve $u$ in the moduli space $\mathcal
M^{J}(\gamma^{t}; \gamma_{i}^{t / k})$, where $k \mid t\leq
N_{\delta}$ and $J$ is an almost complex structure tailored to
$((\mathbb R\times S^1\times D^2)^{\ast},
\tilde{\alpha}_{\delta}^{\ast})$. The index formula can be written
in this way, since $H_{2}(S^1\times D^2; \mathbb Z)=0$ and hence
$<c_{1}(\xi), A>=0$ for all $A\in H_{2}(S^1\times D^2, \mathbb Z)$.
We now use Formulas~\ref{cz1}
and get
\begin{align*}
|\gamma^{t}|-|\gamma_{i}^{t/k}|=-2m-(-2\frac{t}{k}-1)=-2(m-\frac{t}{k})+1,
\end{align*}
and $m=\frac{t}{k}$ for $i=1,\dots,n$; $t\leq N_{\delta}$. Hence,
we can rewrite Equation~\ref{index} as
\begin{align}\label{wrongvalue1}
ind(u)=|\gamma^{t}|-|\gamma_{i}^{t/k}|=-2(\frac{t}{k}-\frac{t}{k})+1=1
\end{align}
for $i=1,\dots,n$ and $t\leq N_{\delta}$. Therefore, Theorem~\ref{transcylinder} and Remark~\ref{translationforsutures} imply that
for every $u\in \mathcal M(\gamma^{t},
\gamma_{i}^{t / k})$ the linearization of the Cauchy-Riemann operator is surjective
at $u$; here $k\mid t\leq N_{\delta}$, $J$ is any
almost complex structure tailored to $((S^1\times D^2)^{\ast}, \tilde{\alpha}_{\delta}^{\ast})$ and $i=1,\dots,n$.

Let $(S^1\times D^2, \Gamma_{long}, U(\Gamma_{long}),
\alpha_{\delta}^{long})$ be a sutured contact solid torus obtained
from $([-1,1]\times D, \Gamma, U(\Gamma), \alpha)$ by using
$\varphi_{X_{H}}$ as a gluing map. Recall that we get $(S^1\times
D^2, \tilde{\Gamma}, U(\tilde{\Gamma}), \tilde{\alpha}_{\delta})$
from $([-1,1]\times D, \Gamma, U(\Gamma), \alpha)$ by using
$\varphi_{-k}=R_{-k}\circ \varphi_{X_{H}}$ as a gluing map. We now
note that $(S^1\times D^2, \Gamma_{long}, U(\Gamma_{long}),
\alpha_{\delta}^{long})$ is a universally tight sutured contact
solid torus with $2nk$ parallel longitudinal sutures, $k>1$, and
such that when one cuts it along the meridian disk the sutures on
the disk are boundary-parallel. This follows from the
gluing/classification result for universally tight contact
structures on a sutured solid torus, see Section 2 in~\cite{Honda}
(more precisely, Corollary~2.3, Theorem~2.5 and Corollary~2.6). The
cylindrical contact homology of this sutured contact manifold is
computed in~\cite{Golovko} and is given by
\begin{align}\label{longcase}
HC^{cyl,h}(S^{1}\times D^2,\Gamma_{long},\xi_{long})
\simeq \left \{
\begin{array}{ll}
\mathbb Q^{nk-1}, & \mbox{for}\ h\geq 1;\\
0, & \mbox{otherwise}.
\end{array}
\right.
\end{align}
Here $\xi_{long}=ker\,\alpha_{long}$.

Note that $(S^1\times D^2, \Gamma_{long}, U(\Gamma_{long}), \alpha_{\delta}^{long})$ has $nk$ hyperbolic orbits $\gamma^{long}_{1},\dots,\gamma^{long}_{nk}$ and one elliptic orbit $\gamma^{long}$. Here $\gamma^{long}_{i}$'s correspond to the equally spaced saddle points of $H$  and $\gamma^{long}$ corresponds to the critical point of $H$ at the center of $D(r_{sing})$. In addition, observe that
\begin{align}\label{longhomology}
[\gamma^{long}_{i}]=[\gamma^{long}]=1\in H_{1}(S^1\times D^2; \mathbb Z).
\end{align}
Finally, note that from Lemma~\ref{orbits} and from the construction of $\gamma^{long}$ and $\gamma_{1}^{long},\dots,\gamma_{nk}^{long}$ it follows that
\begin{align}\label{longaction}
\mathcal A(\gamma^{long})>\mathcal A(\gamma_{i}^{long}),\quad \mathcal A(\gamma_{i}^{long})=\mathcal A(\gamma_{j}^{long})
\end{align}
for $i,j=1,\dots,nk$. Hence, Theorem~\ref{welldefconthom},
Remark~\ref{translationforsutures} together with Fact~\ref{energy},
and Formulas~\ref{longcase}, \ref{longhomology} and \ref{longaction}
imply that $\partial (\gamma^{long})^{s}\neq 0$ for $s>0$; otherwise we come to
contradiction to Formula~\ref{longcase} ($\partial (\gamma^{long})^s= 0$
implies that the exponent of $\mathbb Q$ in Formula~\ref{longcase}
must be $nk+1$). In addition, observe that  $<\partial (\gamma^{long})^{s}, (\gamma^{long}_i)^{s}>\neq 0$ for some $i$ and all $s>0$.

We now take an almost complex structure $J^{long}$ tailored to
$((S^1\times D^2)^{\ast}, (\alpha_{\delta}^{long})^{\ast})$ such
that as a map $\xi^{long}\to \xi^{long}$ it is obtained from some
fixed $J^{cyl}:\xi\to \xi$ which is defined on $([-1,1]\times D,
\Gamma, U(\Gamma), \alpha)$ and satisfies the following properties:
\begin{enumerate}
\item $(J^{cyl})^2=-I$, $d\alpha(J^{cyl}\cdot, J^{cyl}\cdot) = d\alpha(\cdot, \cdot)$, $d\alpha(\cdot,J^{cyl}\cdot)>0$;
\item $J^{cyl}|_{\{1\}\times D} = \varphi_{X_{H}}^{\ast}(J^{cyl}|_{\{-1\}\times D})$
and $J^{cyl}$ is $\frac{2\pi}{nk}$-symmetric, i.e., it is
invariant under $\frac{2\pi}{nk}$-rotation with respect to the
center of $D$.
\end{enumerate}
Here $\xi^{long}=ker\,\alpha_{\delta}^{long}$ and $\xi=ker\,\alpha$.
By saying that $J^{long}$ is obtained from $J^{cyl}$ we simply mean that the gluing procedure with $\varphi_{X_{H}}$ applied to $([-1,1]\times D, \Gamma, U(\Gamma), \alpha)$ transforms $J^{cyl}$ to $J^{long}$.
Since $\xi$ is $\frac{2\pi}{nk}$-symmetric on $([-1,1]\times D, \Gamma, U(\Gamma), \alpha)$, we claim that $J^{cyl}$, which satisfies Properties $(1)$ and $(2)$, exists and that Property $(2)$ is not a serious restriction on $J^{cyl}$. The symmetry of $\xi$ follows from the symmetry of $\beta$ and $X_{H}$, and from the construction of $\alpha$. From the symmetry of $J^{long}$ it follows that $<\partial (\gamma^{long})^s, (\gamma^{long}_i)^s>\neq 0$ for all $i=1,\dots,nk$ and $s>0$.

Now we take $\tilde{J}$ on $(S^1\times D^2, \tilde{\Gamma},
U(\tilde{\Gamma}), \tilde{\alpha}_{\delta})$, which is obtained from
the same $J^{cyl}$ defined on $([-1,1]\times D, \Gamma, U(\Gamma),
\alpha)$ by applying the gluing procedure with
$\varphi_{-k}=R_{-k}\circ \varphi_{X_{H}}$ to $([-1,1]\times D,
\Gamma, U(\Gamma), \alpha)$, and possibly modify it near the
boundary of $(S^1\times D^2, \tilde{\Gamma}, U(\tilde{\Gamma}),
\tilde{\alpha}_{\delta})$ (far from the Reeb orbits) so that it
becomes tailored to $((S^1\times D^2)^{\ast},
(\tilde{\alpha}_{\delta})^{\ast})$. Observe that we can assume that
$J^{long}$=$\tilde{J}$. From the symmetry of $J^{cyl}$ and the form
of the gluing maps for $(S^1\times D^2, \tilde{\Gamma},
U(\tilde{\Gamma}), \tilde{\alpha}_{\delta})$ and $(S^1\times D^2,
\Gamma_{long}, U(\Gamma_{long}), \alpha_{\delta}^{long})$ it follows
that every $J^{long}$-holomorphic curve $u$ which contributes to
$<\partial (\gamma^{long})^{ks}, (\gamma^{long}_{i})^{ks}>\neq 0$
can be modified to a $\tilde{J}$-holomorphic curve $\tilde{u}$ from
$\gamma^{ks}$ to $\gamma_i^s$ by modifying (composing) it with the
rotation about the center of a meridian disk, and hence $<\partial
\gamma^{ks}, \gamma_i^s>\neq 0$.

Observe that this choice of almost complex structures is possible, since Theorem~\ref{transcylinder} and Remark~\ref{translationforsutures} imply that we do not need to require almost complex structures to be generic.
Finally, from Formula~\ref{chaincomplex1} it follows that
\begin{align*}
HC^{cyl,h}_{m}(S^{1}\times
D^2,\tilde{\Gamma},\tilde{\alpha}_{\delta})\simeq\left \{
\begin{array}{ll}
\mathbb Q, & \mbox{for} \ h > 0 \ \mbox{and} \ m=2\lfloor h(-\frac{1}{k}+\delta) \rfloor; \\
\mathbb Q^{n-1}, & \mbox{for}\ k \mid h > 0 \ \mbox{and} \ m=-\frac{2h}{k}-1;\\
0, & \mbox{otherwise}
\end{array}
\right.
\end{align*}
for $h\leq N_{\delta}$.

We now note that $\tilde{\xi}=ker\,\tilde{\alpha}_{\delta}$ is independent
of $\delta$. This follows from the gluing/classification result for universally
tight contact structures on a sutured solid torus, see Section 2
in~\cite{Honda} (more precisely, Corollary~2.3, Theorem~2.5 and Corollary~2.6).
Hence, from
Theorem~\ref{welldefconthom} and Remark~\ref{translationforsutures}
it follows that
\begin{align*}
HC^{cyl,h}(S^1\times
D^2,\tilde{\Gamma},\tilde{\xi})=HC^{cyl,h}(S^1\times
D^2,\tilde{\Gamma},\tilde{\alpha}_{\delta})
\end{align*}
for all $h$ and hence for $h\leq N_{\delta}$,  where $\delta$ is a
small positive irrational number,
\begin{align*}
&HC^{cyl,h}(S^1\times
D^2,\tilde{\Gamma},\tilde{\xi}):=\bigoplus\limits_{m}HC^{cyl,h}_{m}(S^1\times
D^2,\tilde{\Gamma},\tilde{\xi}),\quad \mbox{and}\\
&HC^{cyl,h}(S^1\times
D^2,\tilde{\Gamma},\tilde{\alpha}_{\delta}):=\bigoplus\limits_{m}HC^{cyl,h}_{m}(S^1\times
D^2,\tilde{\Gamma},\tilde{\alpha}_{\delta}).
\end{align*}
Now observe that
$N_{\delta}\rightarrow \infty$ when $\delta\rightarrow 0$. In addition, we note that for fixed $n$, $k$ and
two small positive irrational
numbers $\delta_{1}\neq \delta_{2}$, the sets of closed orbits of
$R_{\tilde{\alpha}_{\delta_{1}}}$ and
$R_{\tilde{\alpha}_{\delta_{2}}}$ are the same, and the corresponding
orbits with the same first homology class $h\leq min\{N_{\delta_{1}},
N_{\delta_{2}}\}$ have the same Conley-Zehnder gradings in the
corresponding complexes. Therefore, for every $0<h\in \mathbb Z=
H_{1}(S^1\times D^2; \mathbb Z)$, there exists $\delta$ such that
\begin{align*}
HC^{cyl,h}_{m}(S^1\times
D^2,\tilde{\Gamma},\tilde{\xi})&=HC^{cyl,h}_{m}(S^1\times
D^2,\tilde{\Gamma},\tilde{\alpha}_{\delta})\\
&\simeq\left \{
\begin{array}{ll}
\mathbb Q, & \mbox{for} \ h > 0 \ \mbox{and} \ m=2\lfloor h(-\frac{1}{k}+\delta) \rfloor; \\
\mathbb Q^{n-1}, & \mbox{for}\ k \mid h > 0 \ \mbox{and} \ m=-\frac{2h}{k}-1;\\
0, & \mbox{otherwise}
\end{array}
\right.
\end{align*}
for $h\leq N_{\delta}$ and hence
\begin{align}\label{withmaslandhom}
HC^{cyl,h}_{m}(S^1\times D^2,\tilde{\Gamma},\tilde{\xi})\simeq\left
\{
\begin{array}{ll}
\mathbb Q, & \mbox{for} \ h > 0 \ \mbox{and} \ m=2\lfloor -\frac{h}{k}+\delta_{k} \rfloor; \\
\mathbb Q^{n-1}, & \mbox{for}\ k \mid h > 0 \ \mbox{and} \ m=-\frac{2h}{k}-1;\\
0, & \mbox{otherwise},
\end{array}
\right.
\end{align}
where $0<\delta_{k}\ll \frac{1}{k}$. Finally, Formula~\ref{withmaslandhom} implies that
\begin{align}\label{lastfinansw}
HC^{cyl,h}(S^{1}\times D^2,\Gamma,\xi)\simeq\left \{
\begin{array}{ll}
\mathbb Q, & \mbox{for} \ k \nmid h > 0;\\
\mathbb Q^{n-1}, & \mbox{for}\ k \mid h > 0;\\
0, & \mbox{otherwise}.
\end{array}
\right.
\end{align}
This completes the proof of Theorem~\ref{mainresult} when $l=1$.

For $l>1$, one can use the same observations as in the case when $l=1$  and show that the only non-zero part of the cylindrical contact homology differential is given by
$<\partial \gamma^{t}, \gamma^{t/k}_{i}>\neq 0$ for $k\mid t\leq N_{\delta}$. This will lead to Formula~\ref{lastfinansw} for all $l$ such that $(k,l)=1$, $k>l>0$.

\begin{remark}
Theorem~1.3 from~\cite{Golovko} and Theorem~\ref{mainresult}
provide the formula for the
sutured version of cylindrical contact homology of $(S^1\times
D^2,\Gamma, \xi)$, where $\Gamma$ consists of
$2n$
parallel sutures of arbitrary slope, $\xi$ is a universally tight contact structure and such that if one
cuts along the meridian disk, the sutures on the disk are $\partial$-parallel.
In particular, this gives a complete calculation of the cylindrical
contact homology of $(S^1\times
D^2,\Gamma, \xi)$, where $\Gamma$ consists of
$2$
parallel sutures of arbitrary slope and $\xi$ is a universally tight contact structure (observe that in this situation there are only two isomorphic (but not isotopic) universally tight contact structures, see Section 2 in \cite{Honda}). These are not all the universally tight contact structures on the solid torus, but all of them can be obtained from the $\#\Gamma=2$ case by successively applying the folding operation.
\end{remark}

\section*{Acknowledgements}
The author is deeply grateful to Ko Honda for his guidance, help and
support. He also thanks Dmytro Chebotarov, Oliver Fabert, Paolo
Ghiggini, Jian He, Michael Hutchings and Mark McLean for helpful
suggestions and interest in his work. In addition, the author is
extremely grateful to Andrew Cotton-Clay for his critical comments
on the first version of the paper. Also, the author is grateful to
the referee of an earlier version of this paper for many valuable
comments and suggestions. Finally, the author thanks the
Mathematical Sciences Research Institute and the organizers of the
``Symplectic and Contact Geometry and Topology'' program for their
hospitality.

\end{document}